\documentclass{article}   	
\usepackage{amsmath,afterpage,bm,amssymb,colortbl,graphicx,enumerate}
\usepackage[hyphens]{url}
\usepackage[colorlinks,citecolor=blue,urlcolor=blue]{hyperref}
\pdfoutput=1
\newenvironment{newreferences}
               {\section*{References}
                \raggedright
                \begin{list}{}{\setlength{\itemsep}{0pt}
                               \setlength{\parsep}{0pt}
                               \setlength{\labelwidth}{0pt}
                               \setlength{\leftmargin}{12pt}
                               \setlength{\labelsep}{0pt}}
                \setlength{\itemindent}{-12pt}
               }{\end{list}}
\def\sp{$\spadesuit$}
\def\cl{$\clubsuit$}
\def\di{$\diamondsuit$}
\def\he{$\heartsuit$}
\def\P{\text{P}}
\def\E{\text{E}}
\allowdisplaybreaks

\title{\textbf{Teaching a University Course on \\ the Mathematics of Gambling}}
\author{Stewart N. Ethier\thanks{\,Department of Mathematics, University of Utah. Email: \href{mailto:ethier@math.utah.edu}{ethier@math.utah.edu}.} \ and Fred M. Hoppe\thanks{\,Department of Mathematics and Statistics, McMaster University. Email: \href{mailto:hoppe@mcmaster.ca}{hoppe@mcmaster.ca}.}}
\date{}

\begin{document}
\maketitle

\begin{abstract}
Courses on the mathematics of gambling have been offered by a number of colleges and universities, and for a number of reasons.  In the past 15 years, at least seven potential textbooks for such a course have been published.  In this article we objectively compare these books for their probability content, their gambling content, and their mathematical level, to see which ones might be most suitable, depending on student interests and abilities.  This is not a book review (e.g., none of the books is recommended over others) but rather an essay offering advice about which topics to include in a course on the mathematics of gambling.\medskip

\end{abstract}

\section{Introduction}

There is anecdotal evidence that a number of college and university courses have been offered on the mathematics of gambling.  Early examples include a course offered in 1970 at what is now California State University, Sacramento (Griffin, 1999, p.~1),
and a course first offered in 1974 at Carleton University (Schneider \& Turmel, 1975).

There are several reasons for such a course.  One would be to train the next generation of casino managers, perhaps in a business program, although the potential audience would be limited.  Another would be to offer liberal arts students a nontechnical course illustrating the applicability of  mathematics, perhaps to meet some graduation requirement.  A third reason would be to give a course in discrete probability, a subject that may seem rather dry to some students, in a more entertaining and accessible way.  A fourth reason would be to offer a course on the mathematics of gambling as a subject in its own right, much as one might give a course on the mathematics of finance.  Such a stand-alone course might even have a probability prerequisite.  

Certainly, there should be no question that the ``theory of gambling'' is a subject worthy of study on its own.  Indeed, as Bachelier (1914, p.~6) wrote over a century ago,\footnote{\,Translation from Dubins \& Savage (2014, p.~vi).}
\begin{quote}
It is almost always gambling that enables one to form a fairly clear idea of a manifestation of chance; it is gambling that gave birth to the calculus of probability; it is to gambling that this calculus owes its first faltering utterances and its most recent developments; it is gambling that allows us to conceive of this calculus in the most general way; it is, therefore, gambling that one must strive to understand, but one should understand it in a philosophic sense, free from all vulgar ideas.
\end{quote}
On the other hand, Laplace (1819, pp.~clxix, v) regarded probability as the subject of primary interest, with gambling having a subsidiary role:\footnote{\,Translation from Laplace (1902, pp.~195, 1).}
\begin{quote}
It is remarkable that a science, which commenced with the consideration of games of chance, should be elevated to the rank of the most important subjects of human knowledge. [$\ldots$]  [T]he most important questions of life [$\ldots$] are indeed for the most part only problems of probability.
\end{quote}
Whether gambling provides motivation for the study of probability, or probability provides tools for the study of gambling, the distinction is primarily one of emphasis.

More-recent examples of gambling courses illustrate the diversity of approaches: 
\begin{itemize}

\item McMaster University: ``Probability and Games of Chance.''  (The original title when the course was proposed was ``Probability and Gambling.'')   Offered seven times since 2009, typically in alternate years.  Average enrollment per term: 52 through 2018; unexpectedly, 155 in 2019.  Prereq.: Introductory probability.  Intended for 3rd- and 4th-year students.  Instructor: Fred M. Hoppe.

\item Carleton University: ``Mathematical Analyses of Games of Chance.''  Offered biennially since 2009.  Prereq.: Permission of school (introductory probability and statistics recommended).  Appears to be intended for 3rd- and 4th-year students.  Instructor: Jason Z. Gao.

\item University of Denver: ``Probability Theory: The Mathematics of Gambling.''  Offered annually since the late 1990s.  Usual enrollment:  40 (full capacity).  Prereq.: Basic statistics.  Intended primarily for juniors and seniors in business.  Instructor: Robert C. Hannum until 2018.

\item Drexel University: ``Optimal Strategies for Repeated Games.''  Offered twice in the early 2010s by ECES (Electrical \& Computer Engineering - Systems).  Average enrollment: 20--25.  Prereq.:  Probability and statistics for engineers.  Intended for seniors and graduate students in engineering.  Instructor: Steven Weber.

\item University of California, Santa Cruz: ``Gambling and Gaming.''  Offered twice a year (in the quarter system) since 2009.  Enrollment 120--150 per term.  Prereq.: Pre-calculus.  Intended for freshmen and sophomores from non-quantitative majors, so as to meet a ``Statistical Reasoning'' requirement.  Instructors: Abel Rodr\'iguez, Bruno Mendes, and others.

\item University of Ottawa: ``Probability and Games of Chance! Poker 101.''  Offered annually.  2019 enrollment: 340.  Prereq.: None.  Intended for first-year students as a science or general elective.  The list of topics comprises game theory, probability theory, history, games (including poker), psychology, and gaming today.  Instructor: Pieter Hofstra.

\item Macquarie University: ``Gambling, Sport and Medicine.''  Offered since 1999.  Annual enrollment: About 300.  Prereq.: None.  Intended for first-year students.  This is a statistics course with applications to the three subjects in the title, with an emphasis on gambling.  Instructor:  John S. Croucher and others.

\item Roanoke College: ``Math of Gambling and Games.''  Offered almost annually since 2009, as part of the ``Intellectual Inquiry Curriculum.''  Enrollment: About 16 per term.  Prereq.: High school algebra.  Intended for students from a broad range of majors, including mathematics.  Instructors: David G. Taylor and others.

\item Albion College: ``Mathematics of the Gaming Industry.''  Offered in Spring 2014 and Spring 2017.  Enrollment both times: 8.  Prereq.: Permission of instructor.   Intended for Mathematics and Computer Science majors/minors who turn 21 by spring break, during which a class trip to Las Vegas is held.  Instructor: Mark Bollman.

\item Stanford University: ``Mathematics and Statistics of Gambling.''  Offered in Spring 2018.  Enrollment: 40, including 10 undergraduates and 15 graduate students from Mathematics/Statistics, and 15 from other departments.  Prereq.: Undergraduate probability and undergraduate statistics.  Intended for graduate students.  Instructor: Persi Diaconis.

\item University of Utah: ``Mathematics of Games of Chance.'' Offered in Summer 2005 as an REU (Research Experience for Undergraduates) program.  Enrollment: 10.  Prereq.:  Permission of instructor following a competitive application process.  Intended for undergraduates interested in mathematical research.  Instructor: Stewart N. Ethier.

\end{itemize}
There are undoubtedly other gambling courses that have not come to our attention.

\begin{figure}[ht]
\begin{center}
\includegraphics[width=4.0in]{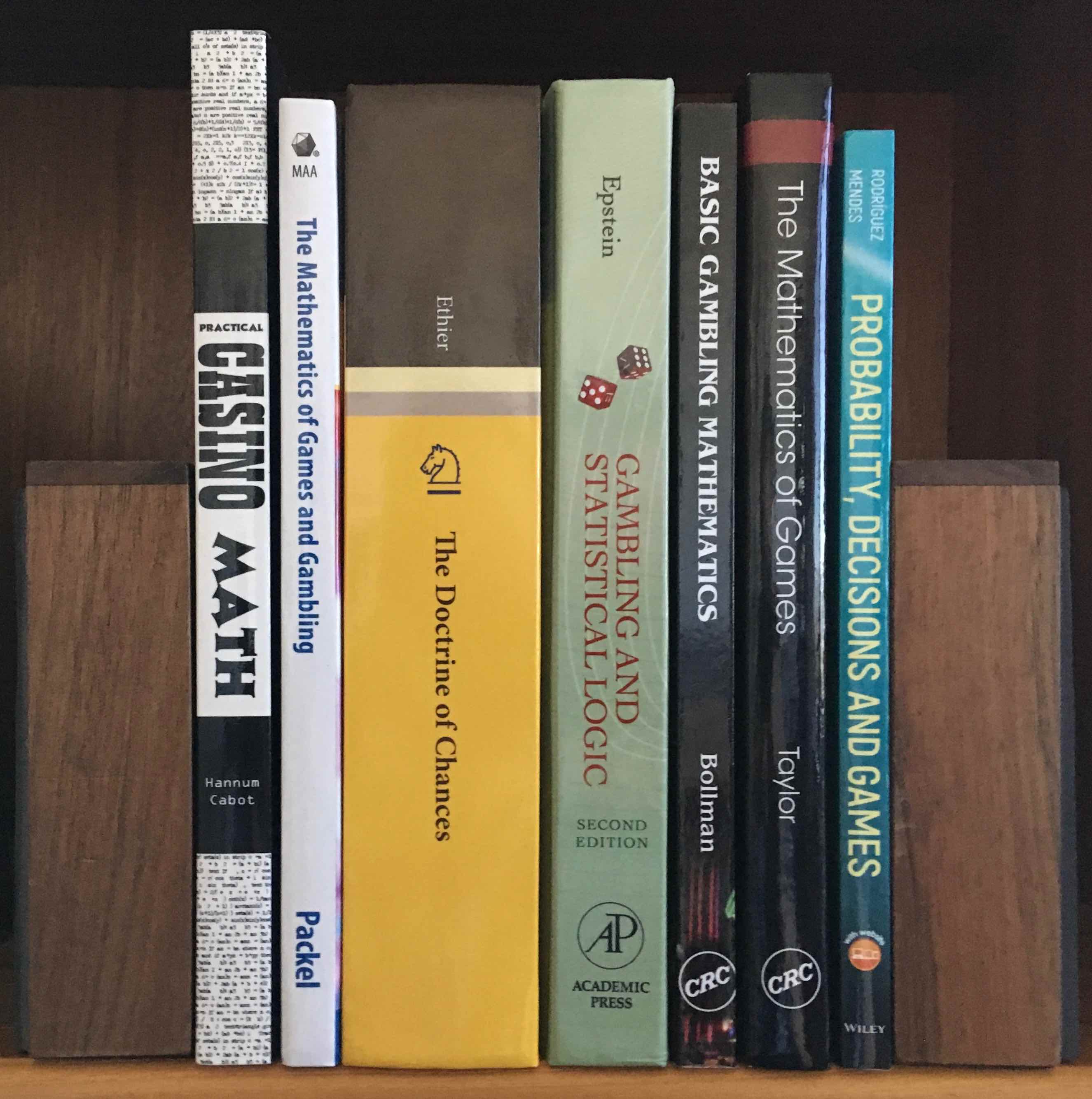}
\vglue5mm
\caption{\label{books}The seven textbooks compared here.} 
\end{center}
\end{figure}

In the past 15 years, at least seven potential textbooks for a course on the mathematics of gambling have been published.  These include the following, listed by year of publication (see also Table~\ref{comparisons} and Figure~\ref{books}):

\begin{enumerate}
\item Hannum, R. C. \& Cabot, A. N. (2005). \textit{Practical Casino Math}, Second Edition. Reno: Institute for the Study of Gambling and Commercial Gaming; Las Vegas: Trace Publications.

\item Packel, E. W. (2006). \textit{The Mathematics of Games and Gambling}, Second Edition.  Washington, D.C.: The Mathematical Association of America.

\item Ethier, S. N. (2010). \textit{The Doctrine of Chances: Probabilistic Aspects of Gambling}.  Berlin--Heidelberg: Springer-Verlag.

\item Epstein, R. A. (2013). \textit{The Theory of Gambling and Statistical Logic}, Second Edition.  Waltham, MA: Academic Press, an imprint of Elsevier.

\item Bollman, M. (2014). \textit{Basic Gambling Mathematics: The Numbers Behind the Neon}.  Boca Raton, FL: CRC Press, an imprint of Taylor \& Francis Group.

\item Taylor, D. G. (2015). \textit{The Mathematics of Games: An Introduction to Probability}. Boca Raton, FL: CRC Press, an imprint of Taylor \& Francis Group.  

\item Rodr\'iguez, A. \& Mendes, B. (2018). \textit{Probability, Decisions and Games: A Gentle Introduction Using R}.  Hoboken, NJ: John Wiley \& Sons. 
\end{enumerate}

\noindent Hannum \& Cabot (2005) is currently out of print (except for a Chinese edition), and Epstein (2013) is currently published only as an e-book.

A few of the books, especially Epstein (2013) and Taylor (2015), have broader coverage than just gambling, but since most textbooks have more material than can be covered in a course, this need not be a liability for a course on the mathematics of gambling.  However, it \textit{is} a liability for Gould's (2016) \textit{Mathematics in Games, Sports, and Gambling}, which has insufficient emphasis on gambling for a course whose focus is gambling, so we exclude it from the list.  Croucher (2003) is excluded because it has no coverage of card games (specifically, baccarat, blackjack, video poker, and poker are not discussed).  B\v arboianu (2013) is excluded because, while its treatment of probability is reliable, its treatment of gambling is not.  Books that are clearly intended for graduate students, specifically Dubins \& Savage (2014) and Maitra \& Sudderth (1996), are also excluded.  We choose to exclude Wilson (1970), Thorp (1984), Haigh (2003),  Bewersdorff (2005), and Shackleford (2019a) because these books are intended for a general audience, and as such are not textbooks, but each is a valuable resource for supplementary material.

\begin{table}[htb]
\caption{\label{comparisons}Comparison of the seven textbooks.}
\catcode`@=\active \def@{\phantom{0}}
\tabcolsep=.14cm
\begin{center}
\begin{tabular}{cccccccc}
\hline
\noalign{\smallskip}
          &  1  &  2  &  3  &  4  &  5  &  6  &  7 \\
author    & Han & Pac & Eth & Eps &  Bol & Tay & Rod \\
edition   & 2nd & 2nd & 1st & 2nd & 1st & 1st & 1st \\
year      & @2005@ & @2006@ & @2010@ & @2013@ & @2014@ & @2015@ & @2018@ \\
\noalign{\smallskip}
\hline
\noalign{\smallskip}
      & xix+ & xiv+ & xiv+ & xiii+ & xi+ & xxi+ & xvii+ \\
pages & 272 & 163  & 744  & 414   & 243 & 288  & 190   \\
      & +27 & +11  & +72  & +37   & +27 & +116 & +25   \\
\noalign{\smallskip}      
chapters & 10 & 7 & 22 & 11 & 7 & 8 & 13 \\
\noalign{\smallskip} 
exercises & 0 & 60 & 468 & 28 & 99 & 110 & 185 \\          
\noalign{\smallskip}                 
answers   & $^1$ & $^2$ & $^3$ & $^4$ & $^5$ & $^6$ & $^7$ \\
\noalign{\smallskip} 
bibliog.\ items  & 58 & 48 & 713 & \phantom{$^*$}663$^*$ & 103 & 43 & 0 \\
\noalign{\smallskip}    
\hline
\noalign{\smallskip}   
\multicolumn{8}{l}{$^1$\,Not applicable.  $^2$\,Answers/hints for selected exercises.  $^3$\,Answers online.}\\
\multicolumn{8}{l}{$^4$\,None.  $^5$\,Answers to odd-numbered exercises. Complete solutions in a}\\
\multicolumn{8}{l}{separate manual.  $^6$\,Answers and selected solutions. $^7$\,Solutions manual for}\\
\multicolumn{8}{l}{odd-numbered problems online.  $^*$\,Some items are counted more than once.}\\
\end{tabular}
\end{center}
\end{table}

In what follows we offer advice about which topics should be included in a course on the mathematics of gambling, and we check the extent to which each the seven textbooks meets our criteria.  For example, video poker is a topic that most students, we believe, would want to learn about, so we make suggestions about which aspects of that subject deserve attention.  We also summarize the video poker coverage of each of the seven textbooks.  In fact, only four of the seven books address video poker at all.

Besides analyzing the coverage of each of these books (in discrete probability and in gambling), we point out any conceptual errors that we have noticed (we found about a dozen) as well as any unconventional approaches.  We emphasize information about the books that can be documented and is not a matter of opinion.  Our aim is to give potential instructors a sense of how appropriate each book would be for the kind of course they may be considering, and the types of students they expect to attract.

We conclude this section with a necessarily brief literature review.  We have found only five papers on this topic, four of which are unpublished (to the best of our knowledge) and from the 1970s, and the fifth of which is more recent.  Each discusses a specific course on the mathematics of gambling.

The Carleton University course first taught by Walter Schneider in 1974 and mentioned above in the opening paragraph was described in some detail in Schneider \& Turmel (1975).  The idea for the course came from Milton Parnes of SUNY Buffalo, who had taught a successful course on gambling in the early 1970s.  Schneider mentioned several problems that he encountered.  One was finding a good textbook.  He eventually settled on Wilson (1970), which is excellent in some areas (e.g., blackjack and baccarat) but is not really a textbook and, from the present-day perspective, is rather dated.  Another problem was that students did not appreciate in-depth coverage of topics such as blackjack, poker, and game theory.  We will explain later why we believe that the students were mistaken.

The next three papers concern a course offered by Larry S. Johnson at Fort Lewis College (Durango, CO) in 1974 (Johnson, 1975, 1977; Gibbs \& Johnson, 1977).  It was a four-week immersion-type course with the second and third weeks spent in Las Vegas.  Textbooks included Wilson (1970) and Thorp (1966a).  The purpose of the course was to motivate students, especially nonmajors, to study mathematics (Johnson, 1975).  There was a wide disparity in the mathematical backgrounds of the students.  This was resolved by putting less emphasis on the mathematics of gambling and more emphasis on non-mathematical aspects of gambling, such as its history.

The fifth paper (Croucher, 2006) describes the motivation for developing the Macquarie University course, ``Gambling, Sport and Medicine,'' mentioned above, which by now has been offered for over 20 years.  The principal aim was to attract more students to the study of statistics, and there is no doubt that this goal was achieved.  To prepare for the course, a textbook was written (Croucher, 2003) and ``props'' were collected, including two working slot machines.  A few years after the course was introduced, it was enhanced using web-based techniques.

\section{Probability content}

Virtually all topics in discrete probability have applicability to gambling, but some may be omitted for the purpose of making the course accessible to a broader audience.  Table~\ref{probability-content} summarizes which of 32 probability topics are included in each of the seven textbooks.

\subsection{Basics}

It does not matter whether probabilities are introduced via axioms, rules, theorems, or definitions.  Ultimately, they are typically evaluated by combining probabilities obtained from the equally-likely-outcomes assumption, using 
\begin{equation}\label{equallylikely}
\P(A)=\frac{|A|}{|S|},
\end{equation}
where $|A|$ is the cardinality of event $A$, and $S$ is the sample space.  This result is so self-evident that Hannum \& Cabot (2005, p.~10) takes it for granted, without formally stating it.  Various notations for the complement of event $A$ include $A^c$ (Ethier, 2010, p.~11), $\overline{A}$ (Epstein, 2013, p.~14; Rodr\'iguez \& Mendes, 2018, p.~2), $A'$ (Bollman, 2014, p.~5), $\text{not }A$ (Packel, 2006, p.~20), and no specific notation (Hannum \& Cabot, 2005, p.~11; Taylor, 2015, p.~11).  Among the basic laws of probability, the only result omitted by most of the textbooks is the inclusion-exclusion law (for $n$ events).  There are gambling applications (e.g., in the game of rencontre, in the derivations of the distribution of dice sums and the distribution of the number of $n$ specified events that occur, and in the analysis of the Fire Bet at craps), but none is essential to the course.

Six of the seven textbooks treat the Chevalier de M\'er\'e's problem of finding both the probability of at least one six in four tosses of a single die and the probability of at least one double six in 24 tosses of a pair of dice (Taylor, 2015, Exercise~1.4, considers only the first part), and four of the seven textbooks treat M\'er\'e's problem of points (Packel, 2006, pp.~14--15; Ethier, 2010, Example~2.1.1; Epstein, 2013, pp.~128--129; Rodr\'iguez \& Mendes, 2018, Section~4.2).  Both problems were subjects of the 17th-century correspondence between Blaise Pascal and Pierre de Fermat, a correspondence that R\'enyi (1972) recreated in a fictional exchange of dated letters between Pascal and Fermat in which Pascal enunciates the laws of discrete probability.

Three of the seven textbooks fail to formalize the basic counting principle, presumably regarding it as self-evident.  Each book covers permutations and combinations, the main distinction being the notation employed.  It is unclear to us why the notation $C_{n,r}$ (Packel, 2006, p.~54) and $_n C_r$ (Bollman, 2014, p.~32) persists, since its only advantage (namely, that it is typewriter friendly) is an anachronism 30+ years after \LaTeX{} was introduced.  Even Hannum \& Cabot (2005, p.~16), which was created in \textit{Word}, uses the standard $\binom{n}{r}$.  For permutations, the notation of choice is $P_{n,r}$ (Hannum \& Cabot, 2005, p.~15; Packel, 2006, p.~54), $_n P_r$ (Bollman, 2014, p.~29; Taylor, 2015, p.~127; Rodr\'iguez \& Mendes, 2018, p.~48), and $P_n^r$ (Epstein, 2013, p.~16).  Ethier (2010, p.~4) departs from this consensus with $(n)_r$, the Pochhammer symbol, though Epstein (2013, pp.~233--234) also uses it occasionally.  One of the few problems that calls for permutations instead of combinations is the enumeration of possible bingo cards.  This is done correctly in Ethier (2010, Problem~14.17), Bollman (2014, Example~3.2.14), and Taylor (2015, Example~8.4), but not in Epstein (2013, p.~157), which uses combinations.

\begin{table}[htb]
\caption{\label{probability-content}Probability content of the seven textbooks.  Entries are page numbers at which the listed topics are addressed.  Page numbers in parentheses indicate that the book treats only a special case.}
\tabcolsep=1mm
\begin{center}
\begin{tabular}{lccccccc}
\hline
\noalign{\smallskip}
      &  1  &  2  & 3 &  4  &  5  &  6  &  7 \\
topic & Han & Pac & Eth & Eps & Bol & Tay & Rod \\
      & 2005 & 2006 & 2010 & 2013 & 2014 & 2015 & 2018 \\
\noalign{\smallskip}
\hline
\noalign{\smallskip}
probability axioms & -- & -- & 12 & 13 & 13 & -- & -- \\
equally likely outcomes & (10) & 17 & 3 & 12 & 15 & 5 & 6 \\
addition law & 12 & 21 & 12 & 13 & 51 & 8 & 11 \\
complementation law & 11 & 21 & 12 & 14 & 14 & 11 & 2 \\
inclusion-exclusion law & -- & -- & 13 & 14 & (51) & (8) & (11) \\
odds & 14 & 18 & 15 & (306) & 46 & 26 & 5 \\
\noalign{\smallskip}
basic counting principle & (70) & (52) & 4 & (16) & 26 & 12 & 8 \\
permutations & 15 & 54 & 4 & 16 & 29 & 127 & 48 \\
combinations & 16 & 54 & 4 & 18 & 33 & 54 & 49 \\
multinomial coefficients & -- & -- & 5 & 20 & -- & 130 & -- \\
\noalign{\smallskip}
independent events & 12 & 21 & 16 & 15 & 54 & 9 & 65 \\
conditional probability & (13) & -- & 17 & 14 & 59 & 92 & 64 \\
multiplication law & 13 & (58) & 18 & 14 & 61 & (93) & 64 \\
total probability law & (110) & -- & 19 & 15 & (69) & -- & 66 \\
Bayes' law & -- & -- & 19 & 15 & -- & -- & 67 \\
\noalign{\smallskip}
random variables & (230) & -- & 21 & 19 & 75 & -- & 15 \\
binomial distribution & 21 & 80 & 22 & 19 & 123 & 115 & 103 \\
hypergeometric distrib. & (82) & (72) & 22 & 24 & 39 & 285 & (50) \\
geometric distribution & -- & -- & 22 & 20 & 146 & 286 & 103 \\
negative binomial dist. & -- & -- & 22 & 179 & -- & 286 & 103 \\
Poisson distribution & 78 & -- & 23 & 20 & -- & -- & -- \\
expectation & 18 & 23 & 27 & 21 & 79 & 32 & 16 \\
variance & 24 & -- & 31 & 24 & 126 & (118) & 21 \\
\noalign{\smallskip}
independent r.v.s & -- & -- & 25 & 22 & -- & -- & 65 \\
multinomial distribution & -- & -- & 26 & 20 & -- & 131 & -- \\
multivariate hypergeom. & -- & (68) & 26 & 19 & -- & (145) & -- \\
conditional expectation & -- & -- & 76 & (178) & -- & (78) & -- \\
expectation of a sum & (26) & -- & 32 & 22 & 81 & -- & 18 \\
covariance & -- & -- & 35 & (278) & -- & -- & -- \\
variance of a sum & (26) & -- & 36 & (57) & -- & -- & (65) \\
\noalign{\smallskip}
law of large numbers & (21) & (88) & 43 & (26) & (17) & -- & 17 \\
central limit theorem & (22) & (85) & 48 & (28) & (128) & (118) & -- \\
\noalign{\smallskip}
\hline
\end{tabular}
\end{center}
\end{table}
\afterpage{\clearpage}

Only three of the textbooks use multinomial coefficients.  Epstein (2013, p.~20) introduces them in connection with the multinomial distribution but never actually uses them as far as we can tell; Taylor (2015, p.~130) introduces them to analyze certain games requiring multiple dice; and Ethier (2010, p.~5) uses them in analyses of video poker and trente et quarante, but also in one elementary application, namely the evaluation of poker-hand probabilities.  

Specifically, the probabilities of one pair, two pairs, three of a kind, full house, and four of a kind need not be treated separately, as is typically done.  Instead, they are all special cases of a single formula, namely
$$
\frac{\binom{13}{d_0,d_1,d_2,d_3,d_4}\binom{4}{1}^{d_1} \binom{4}{2}^{d_2} \binom{4}{3}^{d_3} \binom{4}{4}^{d_4}}{\binom{52}{5}},
$$
where $d_i$ is the number of denominations represented $i$ times in the hand, so that $d_0 +d_1 +d_2 +d_3 +d_4 = 13$ and $d_1 +2d_2 +3d_3 +4d_4 = 5$.  This unified expression seems worth including, regardless of the level of the course.  Perhaps surprisingly, it dates back more than three centuries to Montmort (1708, p.~97; 1713, p.~28); in particular, it predates poker.

\subsection{Conditional probability}

Packel (2006) is unique among the seven textbooks in avoiding the concept of conditional probability.  He formulates the multiplication law as 
\begin{equation}\label{mult-rule-informal}
\P(A\text{ then }B)=\P(A)\P(B),
\end{equation}
\textit{provided that $A$ and $B$ are events in independent successive experiments}, instead of the more general 
\begin{equation} \label{mult-rule}
\P(A\cap B)=\P(A)\P(B\mid A),
\end{equation}
valid for arbitrary events $A$ and $B$.  Consider, for example,
\begin{equation}\label{2aces}
\P(\text{1st card is an ace and 2nd card is an ace})=\frac{4}{52}\,\frac{3}{51}
\end{equation}
(Packel, 2006, p.~17).  It does not follow from \eqref{mult-rule-informal} because the events are not independent, and it does not follow from \eqref{mult-rule} because conditional probability is never defined.  Instead, it follows ``after some thought and consideration of equally likely outcomes'' (p.~16).  Perhaps this means
$$
\P(\text{1st card is an ace and 2nd card is an ace})=\frac{\binom{4}{2}\binom{48}{0}}{\binom{52}{2}}=\frac{4}{52}\,\frac{3}{51},
$$
but Packel (2006, p.~58) eventually reverts to calculation such as \eqref{2aces} without attempting to justify them.

Four of the seven textbooks omit the total probability law, despite its numerous applications in gambling problems.  How then do they evaluate the probability of winning a pass-line bet at craps?  They all effectively use Packel's (2006, Table~3.4) method, namely
\begin{align*}
\P(\text{pass-line win})&=\sum_{k=2}^{12}\P(\text{initial roll is }k\text{ then pass-line win})\\
&=\P(\text{initial roll is }7\text{ or }11)\\
&\qquad{}+\sum_{k=4,5,6,8,9,10}\P(\text{initial roll is }k\text{ then roll }k\text{ before }7),
\end{align*}
which is evaluated by applying \eqref{mult-rule-informal} to each term in the second sum, noting that $\{\text{initial roll}$ $\text{is }k\}$ and $\{\text{roll }k$ $\text{before 7 (starting with the second roll)}\}$ are independent.  Bollman (2014, Example~3.2.6) and Taylor (2015, Table~2.2) adopt this method because they study craps before they introduce conditional probability.

Bayes' law is covered by only three of the seven textbooks, despite its applications in roulette, craps, blackjack, and poker.  For example, if a blackjack player has $\{6,10\}$ vs.\ a playable dealer 10, what is the distribution of the next card dealt?  Here we are implicitly conditioning on the event that the dealer's downcard is not an ace, thereby requiring the use of Bayes' law (Ethier, 2010, Example~1.2.9).  Whether students should be troubled by issues as subtle as this is a choice for the instructor.

An example that can be addressed with or without conditional probability is the now-classic Monty Hall problem.  Popularized by Marilyn vos Savant in the early 1990s, it actually dates back to Selvin (1975a,b).  It is discussed by Ethier (2010, Problem~1.19), Epstein (2013, pp.~143--144), Taylor (2015, pp.~13--15), and Rodr\'iguez \& Mendes (2018, Section~5.1), with the latter confirming its analysis with a simulation (Monty Hall Monte Carlo).

\subsection{Subjective conditional probability}

The conditional probability of $B$ given $A$ is usually defined as the ratio
\begin{equation} \label{cond-prob}
\P(B\mid A)=\frac{\P(A\cap B)}{\P(A)}.
\end{equation} 
To be used effectively, this definition requires that both $\P(A)$ and $\P(A\cap B)$ be computable.
However, typically $\P(B\mid A)$ is specified outside this definition using an interpretation of what conditional probability is meant to represent.  For instance, the second factor in \eqref{2aces},
\begin{equation*}
\P(\text{2nd card is an ace}\mid\text{1st card is an ace})=\frac{3}{51},
\end{equation*}
is found by regarding the conditional probability as the (unconditional) probability of drawing an ace from the 51-card deck obtained by removing an ace from a standard deck (which ace does not matter, and the equally-likely-outcomes assumption applies).  Alternatively but less naturally, we could use the definition \eqref{cond-prob} to derive
\begin{align*}
&\!\!\!\P(\text{2nd card is an ace}\mid\text{1st card is an ace})\\
&=\frac{\P(\text{1st card is an ace and 2nd card is an ace})}{\P(\text{1st card is an ace})}=\frac{\binom{4}{2}\binom{48}{0}/\binom{52}{2}}{4/52}=\frac{3}{51}.
\end{align*}

Because of its connection to betting, a careful introduction to conditioning in a gambling course allows the opportunity to present conditional probability from the subjective or betting perspective, as championed by Bruno de Finetti (1937).  He said that an event $A$ has probability $p = \P(A)\in(0,1)$ \textit{to an individual} if this individual is indifferent as to the choice of the following two bets:
\begin{itemize}
\item (P): pay $1$; receive  $1/p$ if $A$ occurs, nothing if $A^c$ occurs.
\item (H): receive $1$; pay  $1/p$ if $A$ occurs, nothing if $A^c$ occurs.
\end{itemize}
The first corresponds to a player (P) making a bet of size $1$ that pays $1/p-1$ to 1 if $A$ occurs.  
The second corresponds to the house (H) banking the player's bet.

If the bet size is $b > 0$, then all payoffs scale linearly so that (P) and (H) become
\begin{itemize}
\item (P): pay $b$; receive  $b/p$ if $A$ occurs, nothing if $A^c$ occurs.
\item (H): receive $b$;  pay  $b/p$ if $A$ occurs, nothing if $A^c$ occurs.
\end{itemize}
It is easy to check that, by allowing $b < 0$, (P) interchanges with (H) because paying $b$ (positive or negative) is equivalent to receiving $-b$.
 
To say the player is indifferent means that the player does not favor either (P) or (H). 
Indifference implies that both bets are equivalent to
\begin{itemize}
\item Net gain $b(1/p-1)$ if $A$ occurs or net gain $-b$ if $A^c$ occurs, 
\end{itemize}
where $b$ can be positive or negative. (A negative net gain is a net loss.)

Suppose that the player is considering two events, $A$ and $B$, and can place three bets:  A bet on 
$A$ occurring;  or a bet on $A\cap B$ occurring; or a conditional bet on $B$ occurring given that $A$ also occurs, which means that the bet on $B$ takes place only if $A$ occurs.   Otherwise this conditional bet is called off and no money changes hands.   (An example of a conditional bet is the don't pass bet at craps, which is returned to the player if 12 appears on the come-out roll.)

Suppose the player evaluates his or her probabilities of these three events occurring as $\P(A)$, $\P(A\cap B)$, and $\P(B\mid A)$, respectively, which are the player's ``subjective'' probabilities (all in $(0,1)$) of these events occurring.   Let $b_1, b_2, b_3$ be the player's respective bet sizes. 
Summarizing,
\begin{itemize}
\item Bet 1:
	\begin{itemize}
		\item If $A$ occurs, the player wins $b_1(1/\P(A)-1)$.
		\item If $A^c$ occurs, the player wins $-b_1$.
	\end{itemize}
\item Bet 2:
	\begin{itemize}
		\item If $A\cap B$ occurs, the player wins $b_2 (1/\P(A\cap B)-1)$. 
		\item If $(A\cap B)^c$ occurs, the player wins $-b_2$. 
	\end{itemize}
\item Bet 3: If $A$ occurs, then
	\begin{itemize}
		\item If $B$ also occurs, the player wins $b_3 (1/\P(B\mid A)-1)$.
		\item If $B^c$ also occurs, the player wins $-b_3$.
	\end{itemize}
	 But if $A^c$ occurs, then Bet 3 is called off.
\end{itemize}
De Finetti's remarkable insight was that there must be a relationship among the subjective probabilities $\P(A)$, $\P(A\cap B)$, and $\P(B\mid A)$ to prevent a player from becoming a sure loser or a sure winner, which is called avoiding arbitrage or a Dutch book.  This means that it should not be possible to choose wagers $b_1, b_2, b_3$ in such a way that the player will always lose or always win, no matter the random outcome on which the bets are placed.  He called this principle \textit{coherence}.  

To determine what the player wins in any circumstance, partition the sample space into a disjoint union of three events
$$
D_1 = A^c,\qquad
D_2 = A\cap B^c,\qquad
D_3 = A\cap B.
$$
On $D_1$, the player's winnings are
$$
w_1 =  - b_1 - b_2.
$$
On $D_2$, the player's winnings are
$$
w_2 = b_1\bigg(\frac{1}{\P(A)}-1\bigg) - b_2 - b_3.
$$
On $D_3$, the player's winnings are
$$
w_3 = b_1 \bigg(\frac{1}{\P(A)}-1\bigg) + b_2 \bigg(\frac{1}{\P(A\cap B)}-1\bigg) + b_3 \bigg(\frac{1}{\P(B\mid A)}-1\bigg).
$$
Write these three equations in matrix notation by introducing the matrix 
$$
\mathbf M = \begin{pmatrix}
				 -1 & -1  & 0 \\
				1/\P(A)-1 & -1  & -1 \\
				1/\P(A)-1 &  1/\P(A\cap B)-1  & 1/\P(B\mid A)-1
\end{pmatrix}
$$
and column vectors $\mathbf b = (b_1, b_2, b_3)^\textsf{T}$ and $\mathbf w = (w_1, w_2, w_3)^\textsf{T}$, giving
\begin{equation*}
\mathbf w = \mathbf M \mathbf b.
\end{equation*}

If the determinant satisfies $\det (\mathbf M) \ne 0,$ in which case the inverse $\mathbf M^{-1}$ exists, then for any desired winnings vector $\mathbf w$ 
there is a bets vector $\mathbf b$ that will guarantee $\mathbf w$, namely $\mathbf b = \mathbf M^{-1}\mathbf w$.  In particular, if all components of $\mathbf w$ are positive, then the player is a sure winner no matter what outcome occurs.  If  all components of $\mathbf w$ are negative, then the player is a sure loser (and the house a sure winner) no matter what outcome occurs.  To avoid this, de Finetti required that
\begin{equation*}
	\det (\mathbf M) = 0.
\end{equation*}
Evaluating the determinant of $\mathbf M$ with row operations (subtract row 1 from row 2 and then subtract row 1 from row 3) we 
get $1/[\P(A)\P(B \mid A)] - 1/\P(A\cap B) = 0$ which reduces to \eqref{cond-prob} (or \eqref{mult-rule}).  Thus, the formulation \eqref{cond-prob} is consistent with the subjective interpretation of the various probabilities.

This approach is not included in any of the seven textbooks, but maybe it should be.  It requires only a first course in linear algebra.

We implicitly used odds in the development above.   Recall that a player (P) makes a bet of size $1$ that pays $1/p-1$ to 1 if $A$ occurs, or equivalently $1-p$ to $p$.   These are both the true odds and the payoff odds, since $p=\P(A)=p/(p+1-p)$.

\subsection{Random variables}

While each textbook covers expectation, which is arguably the most important concept in the subject, Hannum \& Cabot (2005), Packel (2006), and Taylor (2015) avoid the concept of a random variable.  They can still define the expectation of a distribution, and the house advantage of a wager.  (In fact, Packel uses $X$ for expectation, a letter ordinarily reserved for random variables.)  But there are repercussions for this choice, of which we list several.

\begin{itemize}
\item How does one debunk the claim (from a 1959 Cuban magazine) that, because of a defect in the roulette layout, a one-unit bet on black and a one-unit bet on the third column is a winning play?  With random variables, the gambler's profits $X_1$ and $X_2$ from these two bets satisfy $\E[X_1+X_2]=\E[X_1]+\E[X_2]=-1/19-1/19=-2/19$.  Without random variables, one is obliged to find the distribution of the gambler's profit ($X_1+X_2$), which requires counting the number of black numbers in the third column.  While easy in this case ($\E[X_1+X_2]=(1+2)(4/38)+(-1+2)(8/38)+(1-1)(14/38)+(-1-1)(12/38)=-2/19$; Bollman, 2014, Example~7.1.1), other examples exist where such an accounting would be burdensome.

\item How does one find the house advantage of a pass-line bet with free odds?  The first step is to find the gambler's expected loss.  With the availability of random variables, we can define $X_1$ to be the gambler's loss from a one-unit pass-line bet and $X_2$ to be the gambler's loss from the associated $m$-times free  odds bet ($X_2=0$ if the pass-line bet is resolved on the come-out roll).  Then $\E[X_1]=7/495$ and $\E[X_2]=0$, so $\E[X_1+X_2]=\E[X_1]+\E[X_2]=7/495$.  Without the availability of random variables, one must derive the distribution of the gambler's loss ($X_1+X_2$) and show that its expectation is $7/495$.  This is an extra step that makes the derivation quite a bit more complicated.  The final step in the argument is to divide the expected loss by the expected amount bet, namely $1+(2/3)m$.

\item How does one find the mean $np$ of the binomial$(n,p)$ distribution?  One can find it directly using $k\binom{n}{k}=n\binom{n-1}{k-1}$, as does Bollman (2014, Theorem~4.3.2), or one can use probability generating functions and calculus.  But a much simpler argument, suitable even for the most elementary course, is to write the binomial random variable as a sum of $n$ indicator variables, each with mean $p$.  The same issues arise in finding the mean of the hypergeometric distribution. 

\item How does one state the two gems of probability theory, the strong law of large numbers and the central limit theorem?  Do they not require the random variables
$$
\frac{S_n}{n}\quad\text{and}\quad\frac{S_n-n\mu}{\sqrt{n\sigma^2}}\;?
$$
It is possible to state the Bernoulli case of these theorems (i.e., the normal approximation to the binomial distribution) without random variables, but much is lost by restricting to this special case.  Actually, Hannum \& Cabot (2005) goes beyond the Bernoulli case but does so only very informally.  For example, it uses the self-explanatory notation (p.~26)
\begin{align*}
\text{EV}_\text{win}&=\text{unit wager}\times n\times\text{EV}_\text{per unit},\\
\text{SD}_\text{win}&=\text{unit wager}\times \sqrt{n}\times\text{SD}_\text{per unit}.
\end{align*}
We will return to this point later.
\end{itemize}

Variance too can be computed for a distribution without reference to random variables, and it is important in connection with the central limit theorem, and perhaps especially when comparing slot machines or video poker games.  Only two of the seven textbooks omit the concept of variance (or, equivalently, standard deviation).

The most important special univariate distributions are the binomial, the hypergeometric, and the geometric.  Less important are the negative binomial and the Poisson.  The binomial distribution appears explicitly in the solution of M\'er\'e's problem, the analysis of chuck-a-luck/sic bo, one possible solution of the problem of points, and analysis of the Kelly system, among other applications.  The hypergeometric distribution is ubiquitous in analysis of keno and lottery games, for example.  The geometric distribution appears in the St.\ Petersburg paradox, the martingale system, top-to-random shuffles, the length of the shooter's hand at craps, and the study of casino dead-chip programs, among others.  The negative binomial distribution yields another possible solution of the problem of points but is perhaps not as widely used in gambling calculations.  As a rule, the Poisson distribution does not appear explicitly in gambling problems but only as an asymptotic distribution as some parameters converge.  It therefore deserves less attention.

Next, random vectors and joint distributions can be lightly covered, but should not be omitted.  The two most important multivariate distributions for gambling are the multinomial and the multivariate hypergeometric.  The former is used in connection with biased roulette analysis and games with multiple dice (e.g., Taylor's, 2015, pp.~122--142, study of Yahtzee), and the latter in studies of card games and card counting.  For example, the most natural expression for the probability of a natural in single-deck blackjack is arguably
$$
\frac{\binom{4}{1}\binom{16}{1}\binom{32}{0}}{\binom{52}{2}}.
$$
Independence of random variables is easy, once independence of events is covered.  Variances of sums require covariances except in the case of independent random variables, and because sampling without replacement plays such an important role in card games, one should not assume independence.  

Finally, conditional expectation is a topic that is omitted by most of the seven textbooks, perhaps because a conditional expectation can often be regarded as an unconditional expectation, thereby requiring no new ideas.  Another reason for avoiding this topic is that thinking of a conditional expectation as a random variable is often regarded as an advanced concept, involving $\sigma$-fields and Radon--Nikodym derivatives, more suitable for graduate students.  However, in the setting of discrete random variables, conditional expectations are really quite elementary and nothing to fear.  A careful treatment of conditional expectation yields an expectation version of the total probability law.  This is needed for the century-old result of Brown (1919) evaluating the mean duration of a craps decision.  Such a basic result should not be beyond the scope of a course on the mathematics of gambling.

\subsection{Limit theorems}

The culmination of the probability portion of a course on the mathematics of gambling should include the strong law of large numbers (SLLN) and the central limit theorem (CLT), the two greatest intellectual achievements in the subject.  (Both date back to the 18th century but were not proved under the optimal assumptions until the 20th century.)  But aside from their rich history, these results have important implications for gambling.  The SLLN justifies the frequentist interpretation of probability, and can be used to explain why the house advantage is defined as the ratio of expected loss to expected amount bet, or why the Kelly bettor's fortune grows exponentially.  The CLT quantifies the distribution of deviations between observed and expected results in a precise way.  It is the basis for a casino statistic known as the volatility index.  It also allows approximate evaluation of certain probabilities for which exact computation would be prohibitive.

A nice application of the central limit theorem (for independent but not identically distributed random variables) is given in Hannum (2007), which analyzes the play of an actual online roulette player who had won over a million euros in three weeks.  Was he cheating or was he lucky?

\section{Gambling content}

Many of the simpler gambling games (roulette, craps, keno) can be used to illustrate the concepts of discrete probability, whereas the games that students are typically most interested in are those that offer a better chance of winning (blackjack, video poker, Texas hold'em); these games tend to be more mathematically challenging.  

Textbooks should emphasize contemporary games.  De Moivre (1738) studies bassette (p.~57), pharaon (p.~65), quadrille (p.~83), hazard (p.~135), whisk (p.~147), and piquet (p.~151), but these games are of little interest today, except possibly for historical reasons.  While some of today's casino games have long histories, others are more recent.  The first published reference to Texas hold'em is Livingston (1968).  Video poker is even more recent (late 1970s), and the more popular proprietary games (e.g., Let It Ride, Three Card Poker) date back only to the 1990s.

Table~\ref{gambling-content} summarizes the gambling content of each of the seven textbooks, listing the numbers of pages of coverage of each of 30 topics.  This is not an ideal statistic because some books have more material per page than others, but the table should nevertheless be useful for quick comparisons.

\subsection{Roulette}

One of the simplest casino games is roulette, and it is covered in each of the seven textbooks because it provides a good illustration of the concept of expectation.  Although we focus on the 38-number American wheel, analysis of the 37-number European wheel is nearly identical.  The only important distinction is the partager or en prison rule for even-money bets, which is also present in Atlantic City.  For simplicity, we disregard that rule here.

One useful observation about roulette, noted by only Ethier (2010, p.~461) and Bollman (2014, p.~22), is that the payoff odds of the $m$-number bet ($m=1,2,3,4,6,12,18,24$) are $36/m-1$ to 1.  Thus, if $X_A$ is a random variable denoting the profit from a one-unit bet on a permitted subset $A\subset\{0,00,1,2,3,\ldots,36\}$, then
$$
\E[X_A]=\bigg(\frac{36}{|A|}-1\bigg)\frac{|A|}{38}+(-1)\bigg(1-\frac{|A|}{38}\bigg)=-\frac{1}{19},
$$
so all such bets have the same expectation and house advantage.  We have excluded the notorious 5-number bet whose payoff odds are 6 to 1, rounded down from the $36/5-1=6.2$ to 1 that the formula suggests they should be.

Another observation about roulette is that a one-unit bet on a subset $A\subset\{0,00,1,2,3,\linebreak\ldots,36\}$ (with $|A|$ belonging to $\{1,2,3,4,6,12,18,24\}$) is equivalent to $|A|$ single-number bets of size $1/|A|$ on each of the numbers in $A$ (ignoring the possibility that $1/|A|$ may not be a legal bet size).  To justify this, write $X_A=(36/|A|)I_A-1$, where $I_A$ denotes the indicator of the event that a number in subset $A$ occurs on the next spin.  Then
$$
X_A=\frac{36}{|A|}\,I_A-1=\frac{1}{|A|}\sum_{j\in A}\bigg(\frac{36}{1}\,I_{\{j\}}-1\bigg)=\frac{1}{|A|}\sum_{j\in A}X_{\{j\}}.
$$
This has some useful applications (e.g., Ethier, 2010, p.~466), but here is a trivial one.  Suppose you want to make a bet of $b$ units on $\{0,00,1,2,3\}$ (the 5-number bet).  If $b/5$ is a legal bet size, then instead make five single-number bets of size $b/5$ on 0, 00, 1, 2, and 3.  You will then effectively be paid at the proper odds of 6.2 to 1 if one of those numbers occurs on the next spin.

\begin{table}[htb]
\caption{\label{gambling-content}Gambling content of the seven textbooks.  Entries are the numbers of pages devoted to the listed topics.  Numbers are rounded to the nearest integer.  For example, an entry of \emph{0} means less than half a page of coverage. An en-dash means no coverage.  (Some material may correspond to more than one topic, e.g., poker and matrix/bimatrix games.)}
\tabcolsep=.1cm
\begin{center}
\begin{tabular}{lccccccc}
\hline
\noalign{\smallskip}
      &  1  &  2  & 3 &  4  &  5  &  6  &  7  \\
topic & Han & Pac & Eth & Eps & Bol & Tay & Rod \\
      &  2005 & 2006 & 2010 & 2013 & 2014 & 2015 & 2018 \\
\noalign{\smallskip}
\hline
\noalign{\smallskip}
roulette & 11 & 4 & 28 & 10 & 25 & 9 & 14 \\
craps & 14 & 7 & 46 & 12 & 19 & 12 & 13  \\
keno & 3 & 3 & 20 & 2 & 9 & -- & 0 \\
lotteries & 1 & 10 & 0 & 4 & 8 & 7 & 9 \\
baccarat & 16 & -- & 25 & 2 & 6 & 3 & -- \\
chemin de fer & -- & -- & 11 & 0 & 3 & -- & -- \\
blackjack & 19 & 4 & 49 & 32 & 42 & 37 & 12  \\
video poker & 5 & 7 & 29 & -- & 8 & -- &  -- \\
poker, general & 6 & 7 & 32 & 6 & 6 & 22 & 14  \\
poker, hold'em & 4 & 6 & 28 & 3 & 3 & 2 & 4  \\
\noalign{\smallskip}
slot machines & 22 & 0 & 34 & 1 & 6 & -- & -- \\
bingo & 1 & 0 & 1 & 3 & 4 & 7 & -- \\
big six wheel & 1 & -- & -- & -- & 1 & -- & -- \\
chuck-a-luck/sic bo & 2 & 1 & 0 & 0 & 6 & 1 & -- \\
backgammon & 0 & 13 & 0 & 1 & -- & 5 & -- \\
Casino War & 4 & -- & 2 & -- & 3 & -- & -- \\
Let It Ride & 3 & -- & 11 & -- & 5 & -- & -- \\
Three Card Poker & 3 & -- & 9 & -- & 2 & -- & -- \\
Caribbean Stud & 4 & -- & 2 & -- & 4 & -- & -- \\
pai gow poker & 2 & -- & -- & -- & 4 & -- & -- \\
trente et quarante & -- & -- & 22 & 1 & -- & -- & -- \\
sports betting & 12 & 1 & -- & -- & 8 & -- & -- \\
pari-mutuel betting & 1 & 7 & -- & 9 & 2 & -- & 0 \\
\noalign{\smallskip}
matrix/bimatrix games & -- & 26 & 51 & 27 & -- & 8 & 44 \\
house advantage & 15 & 0 & 36 & -- & 5 & -- & -- \\
betting systems & 3 & 4 & 51 & 6 & 19 & 14 & 7 \\
gambler's ruin & 2 & 6 & 35 & 4 & -- & 9 & -- \\
Kelly system & -- & -- & 36 & 2 & -- & 7 & -- \\
bold play & -- & -- & 40 & 0 & 3 & -- & -- \\
shuffling & 7 & -- & 12 & 9 & 1 & -- & -- \\
\noalign{\smallskip}
\hline
\end{tabular}
\end{center}
\end{table}
\afterpage{\clearpage}

Four of the textbooks consider biased roulette.  An interesting statistical issue is whether the most frequently occurring number in $n$ spins offers a favorable bet.  Ethier (2010, Section~13.2) gives a criterion:  Yes, if the frequency of the most frequent number in $n$ spins exceeds $n/36+c\sqrt{n}$, where $c$ depends on the desired significance level $\alpha$.  Epstein (2013, p.~150) suggests $c=0.48$ if $\alpha=0.05$ and $c=0.40$ if $\alpha=0.20$, though without explanation; these figures are close to those of Ethier (approx.\ 0.49 and 0.41, resp.).  Bollman (2014, Example~4.3.8) suggests the critical value $n/38+0.48\sqrt{n}$ (three standard deviations above the mean), which is the same as Epstein's for $\alpha=0.05$, except for $n/38$ in place of $n/36$.  Rodr\'iguez \& Mendes (2018, Section 7.4) evaluates the probability that a specified number appears at least 270 times in 10,000 spins, but avoids the issue of whether that number was chosen before or after collecting the data.

\subsection{Craps}

Each of the seven textbooks covers craps, at least to the extent of analyzing the pass-line bet.  One step in this analysis to to evaluate the probability of rolling a total of $j$ before a total of 7 in repeated rolls of a pair of dice.  With $\pi_j:=(6-|k-7|)/36$ being the probability of rolling a total of $j$ ($j=2,3,4,\ldots,12$), five of the textbooks give the required probability as $\pi_j/(\pi_j+\pi_7)$, effectively ignoring any result that is not $j$ or 7.  While this is the correct answer, it is somehow less satisfying than the infinite series solution (Ethier, 2010, Example~1.2.1;  Rodr\'iguez \& Mendes, 2018, pp.~79--82),
\begin{align*}
\P(j\text{ before }7)&=\sum_{n=1}^\infty \P(j\text{ before 7, and in exactly }n\text{ rolls})\\
&=\sum_{n=1}^\infty (1-\pi_j-\pi_7)^{n-1}\pi_j=\frac{\pi_j}{\pi_j+\pi_7},
\end{align*}
which requires knowledge of the sum of a geometric series.

The free odds bets at craps deserve attention, being almost unique among casino wagers as conditionally fair bets.  We have already explained that the house advantage of the pass-line bet with free odds is the gambler's expected loss divided by the expected amount bet.  Only Rodr\'iguez \& Mendes (2018) does not cover free odds bets, and Bollman (2014, Example~5.2.3) divides by maximum amount bet instead of expected amount bet.

Craps is usually regarded as one of the simpler casino games, but that does not mean that the mathematics of craps is necessarily trivial.  For example, analysis of the Fire Bet, which pays off according to the number of distinct points made during the shooter's hand, is a challenge.  Ethier (2010, pp.~511--512) gives an explicit formula and Bollman (2014, Example~5.2.5) cites a recursive solution by the Wizard of Odds (Shackleford, 2019b).  Finally, for perhaps the most complicated craps problem ever published, one can refer to a paper by the present authors (Ethier \& Hoppe, 2010), which gives an explicit formula in closed form for the distribution of the length of the shooter's hand at craps.  That paper was motivated by a 2009 incident at the Borgata Hotel Casino \& Spa in Atlantic City in which one Patricia DeMauro rolled the dice 154 times before sevening out, a 5.59 billion to 1 shot.  (This calculation is trickier than it first appears, and a story in \textit{Time} got it wrong, owing to a misunderstanding of the problem; see Suddath, 2009.)  A derivation of the closed formula is probably beyond the scope of any of the seven textbooks.  However, simpler calculations, such as the mean of this distribution (Ethier, 2010, Example~3.2.6, pp.~505--506, p.~510, Problem~15.10; Epstein, 2013, p.~210) or the median (Ethier, 2010, Example~4.1.6; Epstein, 2013, p.~211), are quite manageable, but the variance (Ethier, 2010, p.~511) is complicated.

\subsection{Keno or Lotteries}

Every course should cover keno or lotteries but not necessarily both.  Lotteries have the advantage of life-changing prizes, and live keno has the disadvantage that it is a dying game, not unlike faro a century ago (this does not apply to video keno).  On the other hand, keno, with its way tickets, may be slightly more interesting from the mathematical point of view.  In either case, the key probability distribution is the hypergeometric.  

In this connection, there is one point that requires caution.  The probability of exactly six catches on a 10-spot keno ticket is given by the hypergeometric probability
\begin{equation}\label{keno1}
\frac{\binom{10}{6}\binom{70}{14}}{\binom{80}{20}}=\frac{40{,}583{,}288{,}950{,}923{,}600}{3{,}535{,}316{,}142{,}212{,}174{,}320}=\frac{24{,}869{,}385}{2{,}166{,}436{,}987}\approx0.0114794.
\end{equation}
Some authors (Hannum \& Cabot, 2005, p.~82; Packel, 2006, p.~72; Epstein, 2013, p.~160) prefer to express this as
\begin{equation}\label{keno2}
\frac{\binom{20}{6}\binom{60}{4}}{\binom{80}{10}}=\frac{18{,}900{,}732{,}600}{1{,}646{,}492{,}110{,}120}=\frac{24{,}869{,}385}{2{,}166{,}436{,}987}\approx0.0114794
\end{equation}
because the numerator and denominator are smaller and so easier to evaluate.  But since the casino's choice of 20 numbers is random and the player's choice of 10 numbers need not be, \eqref{keno1} is easy to justify, whereas \eqref{keno2} is not, at least not directly.  Packel's (2006, p.~73) attempted derivation of \eqref{keno2} illustrates the problem:
\begin{quote}
To see how products of combinations arise, consider the case of marking exactly 6 winning numbers [on a 10-spot ticket]. Assume, for the purposes of our reasoning process, that the 20 random numbers have been determined but remain unknown to the player. In how many ways can exactly 6 out of 10 marked numbers appear among the 20 numbers drawn? (In actuality, the player marks his ticket before the numbers are drawn, but our reinterpretation will not affect the results.) There are $\binom{20}{6}$ combinations of 6 marked numbers appearing among the 20 numbers drawn, but for \textit{each such combination} there remain 4 marked numbers to be chosen among the 60 \textit{undrawn} numbers (possible in $\binom{60}{4}$ ways).  Thus there is a total of $\binom{20}{6}\times\binom{60}{4}$ equally likely ways to have exactly 6 of the 10 marked numbers drawn.
\end{quote}
Notice that the basic result \eqref{equallylikely} requires that all outcomes in the sample space be equally likely, not just those corresponding to the event in question.  So this derivation is dubious, especially if the player's 10 numbers are not chosen randomly.  However, it is possible to justify \eqref{keno2} directly (Ethier, 2010, p.~484), that is, without reference to \eqref{keno1}.

We mentioned way tickets as an interesting aspect of keno not found in lotteries.  An especially nice class of way tickets arises as follows.  Let $A_1,A_2,\ldots,A_r$ be mutually exclusive subsets of $\{1,2,\ldots,80\}$, each of size $s\le15$ (so $rs\le80$), and let $t\in\{1,2,\ldots,r\}$ satisfy $st\le15$.  Consider a ticket that bets one unit on $\bigcup_{i\in I}A_i$ for each $I\subset\{1,2,\ldots,r\}$ with $|I|=t$.  This is what might be called an $\binom{r}{t}$-way $st$-spot ticket.  Popular choices include $(r,s,t)=(20,4,2)$ (a 190-way 8-spot ticket) and $(r,s,t)=(10,2,5)$ (a 252-way 10-spot ticket).  Expected payout is easy to evaluate, but the distribution of payout is a challenge.

The distinction between \eqref{keno1} and \eqref{keno2} does not come up in lottery analyses if the sets of numbers chosen by the player and by the lottery are of the same size.  It does come up, however, when the lottery, but not the player, selects a bonus number (see below).  The most interesting feature of a typical lottery is the pari-mutuel payout system employed at least for the higher prize categories.  Specifically, prize pools are shared by the winners, so a wise strategy is to choose a set of numbers that is less likely to be chosen by others.  For example, an arithmetic sequence (e.g., $7,14,21,28,35,42$) would be a poor choice because, if the numbers come up, there will likely be many more winners than usual, resulting is a reduced payout.  A dramatic example of this occurred in Canada's Lotto 6/49 game in 2008, as we explain next.

First, let us state the rules of that game in effect at that time.  The player chooses six numbers from 1 to 49 for a \$2 entry fee (a random choice, called a Quick Pick, is available).  Then Lotto officials randomly draw six numbers as well as a bonus number from 1 to 49.  47\% of sales is dedicated to the Prize Fund. The total amount of \$5 and \$10 prizes are paid from the Prize Fund and the balance of the fund (the Pools Fund) is then allocated between the top four prize categories as indicated in Table~\ref{lotto}. Any amount not won in the top two prize categories is added to the 6/6 Pools Fund for the next draw.  

\begin{table}[htb]
\caption{\label{lotto}Canada's Lotto 6/49 prize structure (PoF = Pools Fund).}\
\catcode`@=\active \def@{\hphantom{0}}
\catcode`#=\active \def#{\hphantom{${,}$}}
\begin{center}
\begin{tabular}{ccccc}
\hline
\noalign{\smallskip}
           &               &                 & 3/19/2008 & 3/19/2008 \\
number of  &  win          & reciprocal of   & winning   & prize per  \\
matches    &               & probability     & tickets   & winner    \\
\noalign{\smallskip}
\hline
\noalign{\smallskip}
6/6$\,^*$     & 80.5\% of PoF & $13{,}983{,}816.0$ & @@#@@0     & $^{\dag}$ \\
5/6 + bonus   & 5.75\% of PoF & @$2{,}330{,}636.0$ & @@#239     & \$1,193.70 \\
5/6, no bonus & 4.75\% of PoF & @@#@$55{,}491.3$   & @@#106     & \$2,223.40 \\
4/6$\,^*$     & 9\% of PoF    & @@#@@$1{,}032.4$   & @$5{,}245$ & @#@\$85.10 \\
3/6$\,^*$     & \$10 prize    & @@#@@@#@$56.7$     & $91{,}935$ & @#@\$10.00 \\
2/6 + bonus   & \$5 prize     & @@#@@@#@$81.2$     & $65{,}843$ & @#@@\$5.00 \\
\noalign{\smallskip}
\hline
\noalign{\smallskip}
\multicolumn{5}{l}{$^{*}$\,Here the bonus number does not play a role.}\\
\multicolumn{5}{l}{$^{\dag}$\,If there had been a single winner, prize would have been \$3,994,188.00.}
\end{tabular}
\end{center}
\end{table}

For example, the probabilities for the top four prize categories are
$$
\frac{\binom{6}{6,0,0}\binom{43}{0,1,42}}{\binom{49}{6,1,42}}, \quad \frac{\binom{6}{5,1,0}\binom{43}{1,0,42}}{\binom{49}{6,1,42}},\quad\frac{\binom{6}{5,0,1}\binom{43}{1,1,41}}{\binom{49}{6,1,42}},\quad\text{and}\quad \frac{\binom{6}{4,1,1}\binom{43}{2,0,41}}{\binom{49}{6,1,42}}+\frac{\binom{6}{4,0,2}\binom{43}{2,1,40}}{\binom{49}{6,1,42}}.
$$
Here we are using a generalization of the multivariate hypergeometric distribution (Ethier, 2010, Eq.~(11.38)).  This is the analogue of \eqref{keno1}.  The analogue of \eqref{keno2} is
$$
\frac{\binom{6}{6}\binom{1}{0}\binom{42}{0}}{\binom{49}{6}},\quad \frac{\binom{6}{5}\binom{1}{1}\binom{42}{0}}{\binom{49}{6}},\quad \frac{\binom{6}{5}\binom{1}{0}\binom{42}{1}}{\binom{49}{6}},\quad\text{and}\quad \frac{\binom{6}{4}\binom{1}{1}\binom{42}{1}}{\binom{49}{6}}+\frac{\binom{6}{4}\binom{1}{0}\binom{42}{2}}{\binom{49}{6}}.
$$
The corresponding probabilities are equal, but the first formulation is easier to justify.

The numbers drawn on March 19, 2008, were $23,40,41,42,44,45$ and the bonus number was 43.  This resulted in a strange anomaly in which the second-category (5/6 + bonus) prize was only about half the size of the third-category  (5/6, no bonus) prize; see Table~\ref{lotto}.  We can estimate that about 6.6 million tickets were sold (and there was no carryover from the previous draw), so the second category would expect to have about three winners, implying that 239 is an extreme outlier.  We speculate that most of the second-category winners chose $40,41,42,43,44,45$, perhaps thinking that no one else would make such a ``nonrandom'' choice.

The idea of choosing numbers that are less likely to be chosen by others is explored in considerable detail in Henze \& Riedwyl (1998).  Among the seven textbooks, Bollman (2014, Section 7.5) gives the most thorough account.

Two other issues are perhaps unique to lotteries.  First, for the typical lottery player, his or her utility function may be nonlinear over the range of possible payouts.  Thus, even if the prize pool has increased to a level that implies a positive expectation, expected utility may nevertheless be negative.  For example, the utility of \$500 million is probably less for most players than 100 times the utility of \$5 million.  The study of lotteries offers a good opportunity to introduce the concept of utility.  A second issue is that big lottery winners are often given the choice of being paid in the form of an annuity or in the form of a reduced lump sum.  It seems that most winners choose the latter, but which choice is mathematically correct?  The answer depends on an actuarial analysis.

\subsection{Baccarat / chemin de fer}

Modern-day baccarat evolved from the 19th century French game chemin de fer, and it is by far the most popular gambling game in the world's largest casino-resort destination, Macau.  Mathematically, chemin de fer is more interesting because it involves some strategy decisions, whereas baccarat does not.  Five of the seven textbooks cover baccarat, but only three treat chemin de fer as well.  The first issue in baccarat is describing Banker's mandatory drawing strategy.  There are several ways to do this.  
\begin{enumerate}[(a)]
\item An $8\times11$ table such as Table~\ref{baccarat-strategy}.  
\item Summary of Table~\ref{baccarat-strategy} row by row as in Table~\ref{bacc-rows}.  
\item Summary of Table~\ref{baccarat-strategy} column by column as in Table~\ref{bacc-cols}. 
\item Summary of Table~\ref{baccarat-strategy} by listing the 13 departures from symmetric play.  (These are the 11 Ss in rows 3--5 and the two Ds in row 6.)
\end{enumerate}  
A case can be made for each approach.  Hannum \& Cabot (2005, p.~100) and Bollman (2014, p.~93) use (b), Epstein (2013, p.~265) and Taylor (2015, p.~266) use (c), and Ethier (2010, p.~598) uses (d).  Approach (b) is apparently the way it is learned by baccarat dealers, the logic being that Banker's two-card hand is seen before Player's third card is dealt.  A potentially useful way to remember Table~\ref{bacc-rows} is, if Banker's two-card total is $x\in\{3,4,5,6\}$ and Player draws a third card $y\in\{0,1,2,\ldots,9\}$, then Banker draws if $2(x-3)\le y\le7$ or if $(x,y)=(3,9)$.

\begin{table}[htb]
\caption{\label{baccarat-strategy}Banker's mandatory drawing strategy at baccarat (D = draw, S = stand).  The shading of D entries is for improved readability.}
\begin{center}
\begin{tabular}{ccccccccccccc}
\hline
\noalign{\smallskip}
Banker's &\multicolumn{11}{c}{Player's third card ($\varnothing$ if Player stands)}\\
two-card &&&&&&&&&&\\
total & 0 & 1 & 2 & 3 & 4 & 5 & 6 & 7 & 8 & 9 & $\varnothing$ \\
\noalign{\smallskip} \hline
\noalign{\smallskip}
0& \cellcolor[gray]{0.85}D & \cellcolor[gray]{0.85}D & \cellcolor[gray]{0.85}D & \cellcolor[gray]{0.85}D & \cellcolor[gray]{0.85}D & \cellcolor[gray]{0.85}D & \cellcolor[gray]{0.85}D & \cellcolor[gray]{0.85}D & \cellcolor[gray]{0.85}D & \cellcolor[gray]{0.85}D & \cellcolor[gray]{0.85}D\\
1& \cellcolor[gray]{0.85}D & \cellcolor[gray]{0.85}D & \cellcolor[gray]{0.85}D & \cellcolor[gray]{0.85}D & \cellcolor[gray]{0.85}D & \cellcolor[gray]{0.85}D & \cellcolor[gray]{0.85}D & \cellcolor[gray]{0.85}D & \cellcolor[gray]{0.85}D & \cellcolor[gray]{0.85}D & \cellcolor[gray]{0.85}D\\
2& \cellcolor[gray]{0.85}D & \cellcolor[gray]{0.85}D & \cellcolor[gray]{0.85}D & \cellcolor[gray]{0.85}D & \cellcolor[gray]{0.85}D & \cellcolor[gray]{0.85}D & \cellcolor[gray]{0.85}D & \cellcolor[gray]{0.85}D & \cellcolor[gray]{0.85}D & \cellcolor[gray]{0.85}D & \cellcolor[gray]{0.85}D\\
3&       \cellcolor[gray]{0.85}D & \cellcolor[gray]{0.85}D & \cellcolor[gray]{0.85}D & \cellcolor[gray]{0.85}D & \cellcolor[gray]{0.85}D & \cellcolor[gray]{0.85}D & \cellcolor[gray]{0.85}D & \cellcolor[gray]{0.85}D & S & \cellcolor[gray]{0.85}D & \cellcolor[gray]{0.85}D\\
4&       S & S & \cellcolor[gray]{0.85}D & \cellcolor[gray]{0.85}D & \cellcolor[gray]{0.85}D & \cellcolor[gray]{0.85}D & \cellcolor[gray]{0.85}D & \cellcolor[gray]{0.85}D & S & S & \cellcolor[gray]{0.85}D\\
5&       S & S & S & S & \cellcolor[gray]{0.85}D & \cellcolor[gray]{0.85}D & \cellcolor[gray]{0.85}D & \cellcolor[gray]{0.85}D & S & S & \cellcolor[gray]{0.85}D\\
6&       S & S & S & S & S & S & \cellcolor[gray]{0.85}D & \cellcolor[gray]{0.85}D & S & S & S\\
7&       S & S & S & S & S & S & S & S & S & S & S\\
\noalign{\smallskip}
\hline
\end{tabular}
\end{center}
\end{table}

The first question a student is likely to have is, Where did Banker's mandatory drawing strategy come from?  To give a partial answer to this question, one could consider a particular case such as $(3,8)$, meaning a Banker two-card total of 3 and a Player third card of 8.  Compute Banker's expectation by drawing and Banker's expectation by standing, and confirm that the latter is larger, perhaps counter-intuitively.  Indeed, the two expectations are 86/1365 and 91/1365, respectively, assuming an infinite-deck shoe, as shown by Wilson (1970, pp.~201--202).  This is consistent with Ethier (2010, Table 5.3).  Bollman (2014, pp.~93--94) computes the standing expectation assuming that the possible Player totals are equally likely, which is true except for 0, and does not compute the drawing expectation.  The origin of Banker's drawing strategy is discussed in some detail by Ethier and Lee (2015b).

\begin{table}[htb]
\caption{\label{bacc-rows}Banker's mandatory drawing strategy at baccarat, as a function of Banker's two-card total.}
\begin{center}
\begin{tabular}{cl}
\hline
\noalign{\smallskip}
Banker's       &   Banker draws $\ldots$ \\
two-card total & \\
\noalign{\smallskip}
\hline
\noalign{\smallskip}
0--2 & always \\
3    & if Player's third card is 0--7, 9, or none \\
4    & if Player's third card is 2--7 or none \\
5    & if Player's third card is 4--7 or none \\
6    & if Player's third card is 6 or 7 \\
7    & never \\
\noalign{\smallskip}
\hline
\end{tabular}
\end{center}
\end{table}

\begin{table}[htb]
\caption{\label{bacc-cols}Banker's mandatory drawing strategy at baccarat, as a function of Player's third card, if any.}
\begin{center}
\begin{tabular}{cc}
\hline
\noalign{\smallskip}
Player's    &   Banker draws   \\
third       &   if Banker's  \\
card        &   two-card total is \\
\noalign{\smallskip}
\hline
\noalign{\smallskip}
0 or 1 & 0--3 \\
2 or 3 & 0--4 \\
4 or 5 & 0--5 \\
6 or 7 & 0--6 \\
8    & 0--2 \\
9    & 0--3 \\
none & 0--5 \\
\noalign{\smallskip}
\hline
\end{tabular}
\end{center}
\end{table}

We strongly encourage the study of chemin de fer, provided the course includes the basics of noncooperative game theory (an issue we will return to).  In the formulation of Kemeny \& Snell (1957), chemin de fer is $2\times 2^{88}$ matrix game.  Player can draw or stand on 5, and Banker's strategy is unrestricted.  Calculations similar to those done for $(3,8)$ above show that Banker should play as in baccarat, except possibly in the cases of $(3,9)$, $(4,1)$, $(5,4)$, and $(6,\varnothing)$.  This reduces the game by strict dominance to $2\times2^4$ and then further to $2\times10$, and the graphical method reveals the $2\times2$ kernel.  The correct solution, in which Player's draw-stand mix when holding 5 is $(9/11,2/11)$ and Banker draws on $(3,9)$, stands on $(4,1)$, draws on $(5,4)$, and uses a draw-stand mix on $(6,\varnothing)$ of $(859/2288,1429/2288)$, appears in Haigh (2003, pp.~217--220) and Bewersdorff (2005, Chap.~42).  Epstein's (2013, p.~265) solution inexplicably contains three errors, despite being stated correctly in his earlier editions, and Bollman (2014, pp.~95--99) studies the game under the assumption that all cards are dealt face up.  In baccarat it does not matter whether cards are dealt face up because drawing rules are mandatory (in practice, cards are dealt face up in mini-baccarat but not in baccarat); but in chemin de fer initial hands must be dealt face down because drawing rules are partly discretionary.

The matrix game solution of Kemeny \& Snell (1957) applies to the \textit{parlor game} chemin de fer, whereas the \textit{casino game} involves a 5\% commission on winning Banker bets, resulting in a bimatrix game.  The Nash equilibrium was found by Ethier \& Lee (2015b) and would be just as suitable for textbook treatment as the Kemeny \& Snell solution.

We recommend doing all chemin de fer calculations under the assumption of an infinite-deck shoe.  The exact $d$-deck analysis, which involves sampling without replacement and composition-dependent two-card hands, is complicated (Ethier \& Gamez, 2013) and likely unsuitable for a textbook.

\subsection{Blackjack}

Each of the seven textbooks covers blackjack but at varying levels of detail.  Let us briefly compare their treatments of basic strategy.  Packel (2006) does not discuss it.  Rodr\'iguez \& Mendes (2018, Section~8.2) gives a portion of basic strategy (ignoring soft hands) and provides motivation, but does not discuss its derivation.  All others provide a complete basic strategy under some set of rules.  Epstein (2013, p.~273--277) does little else, providing no motivation or derivation.  Bollman (2014, Section~6.3) provides some motivation but no derivation.  Hannum \& Cabot (2005, pp.~127--132) provides good motivation but no derivation.  Taylor (2015, pp.~89--99) provides no motivation but does derive basic strategy in a particular case, $\{8,10\}$ vs.\ 10, assuming an infinite-deck shoe (an assumption that avoids the need to condition on the dealer not having a natural).  Ethier (2010, Section 21.2) derives basic strategy in two particular cases, $\{6,10\}$ vs.\ 9 and $\{6,10\}$ vs.\ 10, assuming a single deck, then provides an algorithm for recursively deriving composition-dependent basic strategy, and finally reports the results.

There is no question that deriving basic strategy is complicated, but we do not consider that a good reason to avoid it.  At a minimum, one or two particular cases of basic strategy should be analyzed in some detail to give a sense of the magnitude of the problem.  Assuming an infinite-deck shoe does simplify matters considerably and might permit a more complete analysis.  See Werthamer (2018, Section~7.1) and Bewersdorff (2005, Chap.~17) for viable approaches.  The recursive algorithm mentioned above dates back to Manson et al.\ (1975) and was used by Griffin (1999, p.~172).  In the end it requires a rather elaborate computer program, which raises the question of whether this approach is suitable for a textbook.

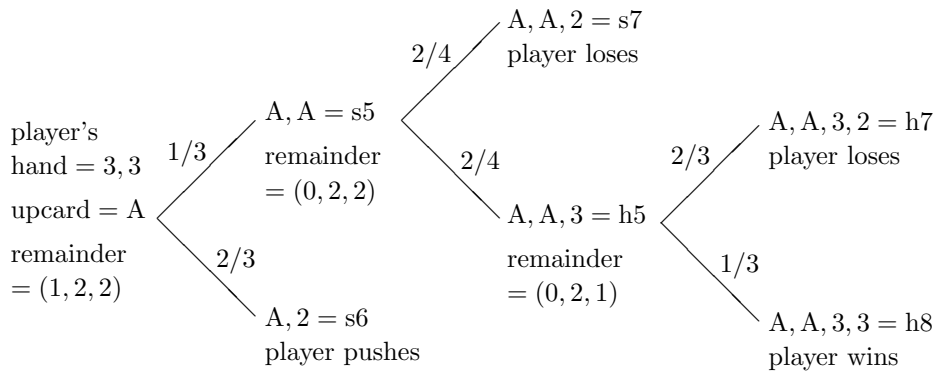
\begin{figure}[htb]
\setlength{\unitlength}{0.65cm}
\begin{picture}(19,8)(.8,0)
\put(1,4.6){player's}
\put(1,3.9){$\text{hand}=3,3$}
\put(1,3){$\text{upcard}=\text{A}$}
\put(1,2.1){$\text{remainder}$}
\put(1,1.4){$=(1,2,2)$}
\put(4,3){\line(1,1){2}}
\put(4,3){\line(1,-1){2}}
\put(4.2,4.2){1/3}
\put(5.2,2){2/3}
\put(6.2,5){$\text{A},\text{A}=\text{s}5$}
\put(6.2,4.1){$\text{remainder}$}
\put(6.2,3.4){$=(0,2,2)$}
\put(6.2,0.8){$\text{A},2=\text{s}6$}
\put(6.2,0.1){$\text{player pushes}$} 
\put(9,5){\line(1,1){2}}
\put(9,5){\line(1,-1){2}}
\put(9.2,6.2){2/4}
\put(10.2,4){2/4}
\put(11.15,6.9){$\text{A},\text{A},2=\text{s}7$}
\put(11.15,6.2){$\text{player loses}$} 
\put(11.15,2.9){$\text{A},\text{A},3=\text{h}5$}
\put(11.15,2){$\text{remainder}$}
\put(11.15,1.3){$=(0,2,1)$}
\put(14.3,2.9){\line(1,1){2}}
\put(14.3,2.9){\line(1,-1){2}}
\put(14.5,4.1){2/3}
\put(15.5,1.9){1/3}
\put(16.5,4.8){$\text{A},\text{A},3,2=\text{h}7$}
\put(16.5,4.1){$\text{player loses}$} 
\put(16.5,0.7){$\text{A},\text{A},3,3=\text{h}8$}
\put(16.5,0.0){$\text{player wins}$} 
\end{picture}
\vglue5mm
\caption{\label{tree}The tree diagram used to evaluate the standing expectation with $\{3,3\}$ vs.\ ${\rm A}$ at single-deck snackjack.  (A remainder of $(1,2,2)$ means that $1$ ace, $2$ deuces, and $2$ treys remain. {\rm s5} and {\rm h5} stand for soft $5$ and hard $5$.)  We find that $\E=-2/9$.  Notice that we are conditioning on the dealer not having a natural (i.e., the dealer's downcard is not a trey).  The corresponding tree diagram for standing with $\{10,10\}$ vs.\ ${\rm A}$ at single-deck blackjack would have $5994$ terminal vertices (Griffin, 1999, p.~158) instead of four.}
\end{figure}

Ethier \& Lee (2019) proposed an alternative approach.  Specifically, they showed that basic strategy can be derived by hand (using the recursive algorithm just mentioned) for a toy model of blackjack called \textit{snackjack} (so-named by Epstein, 2013, p.~291).  The eight-card deck comprises two aces (value 1 or 4), two deuces (value 2), and four treys (value 3).  The target total is 7, not 21, and ace-trey is a natural.  The dealer stands on 6 and 7, including soft totals, and otherwise hits. The player can stand, hit, double, or split, but split pairs receive only one card per paircard (like split aces in blackjack), and there is no insurance.  To see why calculations by hand are feasible, there are only 32 decision points and 17 dealer drawing sequences in single-deck snackjack, compared to 19,620 and 48,532, respectively, in single-deck blackjack.  To see more visually the simplicity of snackjack, refer to Figure~\ref{tree}.  We believe that this toy model has pedagogical value, which could be exploited in a course on the mathematics of gambling.

The next main topic in blackjack is card counting, and there is a lot that can be said.  Topics might include the fundamental theorem of card counting, card counting for insurance decisions, the effects of removal, balanced and unbalanced card-counting systems, betting correlation, strategy variation, the Illustrious 18, and so on.  Among the seven textbooks, Bollman (2014, Section~6.4) has the most thorough coverage.

\subsection{Video poker}

Video poker shares two features with blackjack.  First, these are the only two house-banked casino games in which a skillful player can have a positive expectation under ordinary circumstances (ruling out, for example, games with progressive jackpots).  Second, the mathematical analysis of both games is complicated and computer-intensive.  The first property implies that both games are essential components of any course on the mathematics of gambling, despite the second property.  

A sound approach to video poker is to focus on one specific game, such as the classic 9/6 Jacks or Better.  (Here 9/6 refers to the fact that the payoff odds for a full house are 9 for 1 and those for a flush are 6 for 1.)  See Table~\ref{JoB} for the pay table.  Hannum \& Cabot (2005, pp.~159--160) considers the standard 9/6 game and gives the payoff distribution under optimal play, but does not describe optimal play.  Ethier (2010, Section 17.1) analyzes the standard 9/6 game and gives the payoff distribution under optimal play, as well as exact optimal play.  Bollman (2014, pp.~99--103) studies the 9/6 game, except that a royal flush pays 940 for 1, and gives what is said to be optimal play but is only an approximation to it.  (For example, his hand-rank table ranks ``3-card straight flush'' ahead of ``1 or 2 high cards.''  But A\he-3\di-5\cl-7\cl-9\cl, for example, is correctly played by holding only the ace, not the three-card straight flush.)  Packel (2006, pp.~62--67) treats 8/5 Jacks or Better, except that a royal flush pays 2500 for 1, and optimal play is not described except for determining which hands are pat hands.  

While 9/6 Jacks or Better (optimal expected return 99.5439\%) and 8/5 Jacks or Better (optimal expected return 97.2984\%) are widely available, the two games just mentioned with royal flush bonuses are not currently available (see \url{https://www.vpfree2.com/video-poker/pay-table/jacks-or-better}).  A 2500-for-1 payout on a royal flush is conceivable in a game with a progressive jackpot.  A good reason to choose the standard 9/6 game for analysis is that tutorial software is available online at \url{https://www.freeslots.com/poker.htm}.  This software also has the ability to analyze an arbitrary specified hand.

\begin{table}[htb]
\caption{\label{JoB}Jacks or Better pay tables.  Main entries are the returns from a one-unit (max-coin) bet.}
\catcode`@=\active \def@{\hphantom{0}}
\tabcolsep=3.3mm
\renewcommand{\arraystretch}{1.}
\begin{center}
{\begin{tabular}{lccc}
\hline
\noalign{\smallskip}
rank & 9/6 game\,$^1$               & 9/6 game,$^2$    & 8/5 game,$^3$  \\
     & \phantom{royal bonus0}  & royal bonus  & royal bonus \\
\noalign{\smallskip}\hline
\noalign{\smallskip}
royal flush & 800 & 940 & 2500 \\
straight flush & @50 & @50 & @@50 \\
four of a kind & @25 & @25 & @@25  \\
full house & @@9 & @@9 & @@@8 \\
flush & @@6 & @@6 & @@@5 \\
straight & @@4 & @@4 & @@@4  \\
three of a kind & @@3 & @@3 & @@@3  \\
two pairs & @@2 & @@2 & @@@2  \\
one pair, jacks or better & @@1 & @@1 & @@@1  \\
other & @@0 & @@0 & @@@0  \\
\noalign{\smallskip}\hline
\noalign{\smallskip}
optimal expected return & 99.5439\% & 99.9030\% & 102.3886\% \\
std.\ deviation (opt.\ play) & 4.4175 & 5.1835 & 14.1466 \\
royal probab.\ (opt.\ play) & 1/40,390.5 & 1/38,077.8 & 1/31,746.0 \\
\noalign{\smallskip}\hline
\noalign{\smallskip}
\multicolumn{4}{l}{$^1$\,Studied by Hannum \& Cabot (2005, p.~160) and Ethier (2010, p.~546).}\\
\multicolumn{4}{l}{$^2$\,Studied by Bollman (2014, p.~99).  $^3$\,Studied by Packel (2006, p.~63).} \\
\end{tabular}}
\end{center}
\end{table}
What should a textbook treatment include?  It should certainly discuss equivalence of hands (via permutation of suits) and the fact that there are 134,459 equivalence classes.  It should consider a few example hands and evaluate expected returns under two or more drawing strategies.  It might discuss the concept of penalty cards. It should give the exact (or a nearly) optimal strategy in terms of a hand-rank table, as well as the distribution of return under optimal play.  It might point out that that this strategy is, in effect, uniquely optimal (in 9/6 Jacks or Better but not necessarily in other games).  An extended treatment could give, at least in part, the distribution of the optimal conditional expected return, given the initial hand, a distribution (for 9/6 Jacks or Better) of a random variable with 1153 distinct values (or 387 distinct values if garbage hands are ignored) that was only recently evaluated, apparently for the first time.  For details, see Ethier et al. (2019).

We mention in passing a result from Hannum \& Cabot (2005, p.~162) that could be misinterpreted: 
\begin{quote}
Without delving further into the mathematics, suffice it to say that for multi-play video poker, there is no change in the expected value (house edge) but the variance (volatility) increases.
\end{quote}
The most natural way to compare $n$-play video poker with single-play video poker is to divide a one-unit (max-coin) bet evenly among the $n$ plays, thereby maintaining a bet size of one unit.  Then the mean return is constant in $n$ and the variance is decreasing in $n$.  A less natural comparison is to make a one-unit (max-coin) bet on each play, thereby betting a total of $n$ units.  Then the mean return \textit{per unit bet} is constant in $n$ and the variance is increasing in $n$.  It is this second formulation that was apparently intended.

\subsection{Poker} 

One can study the principles of generic poker, such as hand counting, poker variants (seven-card stud, Omaha, etc.), poker models, pot odds, bluffing, and game theory, but most students will be more interested in learning about Texas hold'em, a game whose popularity has exploded in the past two decades.  Unfortunately, the mathematics of hold'em is either trivial (e.g., the probabilities of making a hand at the turn or the river) or extremely complicated (e.g., calculations of the sort described in the next paragraph).

Hannum et al. (2012) introduced the game of \textit{chance hold'em}, which differs from Texas hold'em in that all players must call the big blind, and then all players must check through to the showdown, regardless of the board.  Their object was to isolate the chance aspect of Texas hold'em, and thereby show, using an online poker database, that skill predominates.  Heads-up chance hold'em can be fully analyzed, but it requires a substantial computing effort (simulation is also feasible but much less reliable).  Ethier (2010, p.~707--715) gives the details for the case of A\sp-K\sp{} vs.\ 8\he-8\di, a calculation that ultimately must be confirmed by computer.  But this is only one of 47,008 heads-up matchups (or equivalence classes of such).  This leads to a ranking of initial hands according to their expected net gain against a random hand in heads-up play, in units of the big blind.  A similar ranking of initial hands in a game with 10 players requires simulation.  See Table~\ref{Texas} for a partial ranking.  Finally, Ethier (2010, pp.~721--725) gives two examples (borrowed from Guerrera, 2007, pp.~62--73) requiring a Bayesian analysis of a call-or-fold decision in a specific game in which the history of the betting is known.  One of them uses the heads-up probabilities of the sort described above.  

The other six textbooks provide only brief coverage of Texas hold'em.

\begin{table}[htb]
\caption{\label{Texas}Ranking of pocket pairs among the 169 distinct initial hands in Texas hold'em.  Columns 2 and 3 rank the hands according to expected net gain vs.\ a random hand in a 2-player (heads-up) game, in units of the big blind; these were determined by exact computation.  Columns 4 and 5 rank the hands according to expected net gain vs.\ random hands in a 10-player game, in units of the big blind; these were determined by a simulation of more than 20.5 billion games (Shackleford, 2016).  In both cases it was assumed that all players call the big blind and then check through to the showdown, regardless of the board.}
\catcode`@=\active \def@{\hphantom{0}}
\tabcolsep=2mm
\renewcommand{\arraystretch}{1.}
\begin{center}
\begin{tabular}{ccccccc}
   && \multicolumn{2}{c}{2-player game@@} && \multicolumn{2}{c}{10-player game@@} \\
\noalign{\smallskip}
\hline
\noalign{\smallskip}
pocket &@@&    & expectation vs. &@@&    & expectation vs. \\
pair && rank & random hand     && rank & random hands     \\                    
\noalign{\smallskip}\hline
\noalign{\smallskip}
AA && 1 & 0.704074  && 1 & 2.1071 \\
KK && 2 & 0.647914  && 2 & 1.6079 \\
QQ && 3 & 0.598503  && 3 & 1.2224 \\
JJ && 4 & 0.549389  && 5 & 0.9318 \\
TT && 5 & 0.500236  && 13 & 0.7145 \\
99 && 6 & 0.441145  && 17 & 0.5559 \\
88 && 7 & 0.383261  && 26 & 0.4471 \\
77 && 9 & 0.324720  && 35 & 0.3647 \\
66 && 17 & 0.265695 && 46 & 0.3017 \\
55 && 27 & 0.206498 && 52 & 0.2379 \\
44 && 48 & 0.140456 && 56 & 0.2155 \\
33 && 66 & 0.073862 && 60 & 0.2002 \\
22 && 87 & 0.006680 && 62 & 0.1933 \\
\noalign{\smallskip}
\hline
\end{tabular}
\end{center}
\end{table}

\subsection{Other games}

We are going to take the position, which is certainly arguable, that all other gambling games should be ignored.  There is simply too much material in the seven topics just described (roulette, craps, keno/lotteries, baccarat/chemin de fer, blackjack, video poker, and poker) to spend time with lesser games.  That's not to say that Let It Ride, for example, lacks mathematical interest.  Or that slot machines or sports betting lacks popularity.  We are simply recognizing that thorough coverage of the principal topics will have greater impact on students than spotty coverage of many topics.  

As for slot machines, the decision to ignore them is actually a consequence of a feature they uniquely enjoy:  The information needed to analyze slots  is highly classified and in particular not available to textbook authors or university instructors.  Ethier (2010, pp.~436--441) analyzes a par sheet made available by Bally for an outdated machine.  The same author also acquired par sheets for some modern five-reel video slots through a confidential source, but his publisher (Springer) was unwilling to risk a lawsuit by publishing that information.  As long as the gambling public continues to patronize games that deny them knowledge of the edge they are up against, the practice of restricting access to that information will persist.

As for sports betting and pari-mutuel betting, analysis of these games is more empirical/statistical than probabilistic, which is why we prefer the seven games mentioned above.

\subsection{Game theory}

We believe that a course on the mathematics of gambling should include the basics of game theory, not because that is a beautiful subject, which it certainly is, but because it is needed to analyze the games of chemin de fer and poker.  Other gambling games that require game theory include le her, an 18th-century card game that was the very first game analyzed using game theory, more than two centuries before game theory existed as a subject; super pan 9, a baccarat-like card game; and baccarat banque, a three-person zero-sum game (Ethier \& Lee, 2015a).  As explained earlier, chemin de fer is a $2\times 2^{88}$ matrix (or bimatrix) game that is reduced to $2\times2^4$, then to $2\times10$, by strict dominance and then solved graphically to reveal the mixed-strategy solution.  Le her is easier, being effectively a $14\times14$ matrix game that is reduced to $2\times2$ by strict dominance, resulting in a mixed-strategy solution.  Super pan 9 is also easier than chemin de fer (but less interesting), being initially a $2\times2^{20}$ matrix game that is reduced to $2\times2$ by strict dominance, then solved by noticing the presence of a saddle point.  We have already discussed poker endgame.  Further examples of poker models to which game theory applies can be found in Chen \& Ankenman (2006).
 
Four of the seven textbooks have extensive coverage of matrix and/or bimatrix games (the subset of game theory most relevant to gambling).  Packel (2006, Chap.~6) treats matrix and bimatrix games, cooperative $n$-person games, and sequential games of perfect information, all in a chapter titled ``Elementary Game Theory,'' and then applies the theory to a poker endgame situation.  Rodr\'iguez \& Mendes (2018, Chaps.~10--13) covers the same topics as Packel except for cooperative $n$-person games.  Ethier (2010, Chap.~5) restricts attention to matrix games, while Epstein (2013, pp.~30--37 and Chap.~10) treats matrix games as well as providing extensive coverage of sequential games of perfect information (e.g., tic-tac-toe, nim, checkers, chess, go).

Restriction to matrix games (rather than bimatrix games) is undesirable because chemin de fer and poker are both examples of matrix games that are more properly regarded as bimatrix games.  We have already mentioned that, with a 5\% commission on Banker wins, the casino game chemin de fer is a bimatrix game, whereas, with no commission, the parlor game chemin de fer is a matrix game.  Poker endgame is a constant-sum bimatrix game, which can be treated as a matrix game if necessary, but it would be clearer not to.  Let us clarify this point.

Basic endgame is described in Ethier (2010, p.~694) (essentially from Ferguson \& Ferguson, 2007) as follows:
\begin{quote}
Two players, 1 and 2, ante $a$ units each, where $a>0$ is specified.  Player 1 then draws a card that gives him a winning hand with probability $P$, where $0<P<1$, and a losing hand otherwise.  Both players know the value of $P$ but only player 1 knows whether he has a winning hand.  It is assumed that play begins with a check by player 2.  Player 1 may then check or bet $b$ units, where $b>0$ is specified.  If he checks, there is a showdown.  If he bets, player 2 may fold or call.  If player 2 folds, player 1 wins the pot.  If player 2 calls, there is a showdown.
\end{quote}
It turns out that player 1 should never check with a winning hand (such strategies are dominated), so the payoff matrix has the following form, with player 1 being the row player and player 2 being the column player:
\begin{equation*}
\catcode`@=\active \def@{\phantom{0}}
\begin{array}{ccc}
& \multicolumn{2}{c}{
\begin{array}{lr}
@@@{\rm fold\ if} &  {\rm call\ if}@@@ \\
{\rm player\ 1\ bets}@ & {\rm player\ 1\ bets}\\
\end{array}
}\\
\noalign{\smallskip}
\begin{array}{l}
{\rm check\ if\ loser{,}\ bet\ if\ winner} \\
{\rm bet\ if\ loser{,}\ bet\ if\ winner} \\
\end{array}
&
\multicolumn{2}{c}{
\left(
\begin{array}{lr}
@@(a_{11},b_{11})@@ & @@(a_{12},b_{12})@@ \\
@@(a_{21},b_{21})@@ & @@(a_{22},b_{22})@@ \\
\end{array}
\right)
}
\end{array}.
\end{equation*}
The $a_{21}$ entry, for example, is the payoff to player 1 if player 1 bets regardless of his status as winner or loser and player 2 folds.  It is therefore the value of the pot, $2a$, since the pot is regarded as belonging to neither player.  The $b_{21}$ entry is the payoff to player 2 under the same circumstances, which is 0, and therefore this is a general-sum (in fact, constant-sum) game.  To study this game using the theory of matrix games, Ethier (2010, pp.~604--607) simply redefined $b_{ij}:=-a_{ij}$, maintaining the correct payoff for player 1 but reducing player 2's payoff by $2a$.  This gives the right answers but is a little contrived.  An alternative, and perhaps more natural, approach is to regard the money in the pot as belonging to the player that contributed it, in which case $a_{21}=a$ and $b_{21}=-a$ (Ferguson \& Ferguson, 2007).  The most natural approach, we believe, is to formulate the game as a bimatrix game.

\subsection{House advantage}

The house advantage of a wager is the most important casino statistic.  It can be expressed as either a fraction or a percentage.  The definition varies slightly from one source to another because there are several aspects of it that are controversial.  Here we mention three issues.  
\begin{enumerate}
\item Should pushes be included or excluded?  For example, should the house advantage of the don't pass bet at craps be $27/1980\approx0.0136364$ (pushes included) or $27/1925\approx0.0140260$ (pushes excluded)?  
\item In compound wagers should the house advantage be computed with respect to the initial amount bet or the total amount bet?  For example, for the ante-play wager at Three Card Poker, is the house advantage (pushes included)
$$
\frac{686{,}689}{20{,}358{,}520}\approx0.0337298\quad\text{or}\quad\frac{686{,}689}{34{,}084{,}400}\approx0.0201467,
$$
that is, expected loss per initial amount bet or expected loss per expected total amount bet?  
\item Should amount bet be the amount ``in play'' or the amount placed at risk?  For example, assuming that a \$5 20-spot keno ticket returns at least \$1, regardless of what happens, should the bet size be \$5 (amount ``in play'') or \$4 (amount placed at risk)?  
\end{enumerate}
These issues are addressed in some detail by Hannum \& Cabot (2005, pp.~56--58) and Ethier (2010, pp.~236--238). Incidentally, our preferred answers are as follows: 1. It depends on the game (e.g., pushes included in baccarat, blackjack, and video poker, pushes excluded in craps and faro). 2. Expected total amount bet, except in blackjack.  3. Amount placed at risk.

Another interesting issue is how the house advantage is affected by a dead-chip (or nonnegotiable-chip) program; see Hannum \& Cabot (2005, pp.~216--231) for a thorough discussion.

\subsection{Betting systems}

Betting systems is a topic worthy of some attention.  Here we are referring to the classic systems for repeated even-money wagers that go by names such as the martingale, Fibonacci, Labouchere, and d'Alembert.   In the martingale system, one doubles one's bet after each loss until finally achieving a win.  Its analysis is straightforward.  

The Fibonacci and Labouchere systems are closely related and mathematically interesting.  In both systems, the bettor keeps an ordered list of numbers on a scoresheet.  The next bet is always the sum of the first and last numbers on the list in the Labouchere system, or the sum of the last two numbers on the list in the Fibonacci system.  (An obvious exception is made if the list has only one number.) After a win, the first and last numbers are canceled (Labouchere) or the last two numbers are canceled (Fibonacci), whereas after a loss, the amount just lost is appended to the list as a new last term.  When the list becomes empty, the bettor stops.  The reason for the name Fibonacci is that, if the initial list is, for example, $1,1,2,3$, subsequent terms and bet sizes are always Fibonacci numbers.  An important distinction is that while the Fibonacci system is amenable to analysis, the Labouchere system is not because of its history dependence.  Epstein (2013, pp.~52--57) gives good coverage, but accuracy is uneven.  For example, the calculation on p.~54 purporting to show the size of the bettor's bankroll needed for a 0.99 probability of success at the Labouchere system at a fair game starting with an initial list of $1,2,3$ is erroneous, and we do not know how to fix it; this material is new to the second edition.  Bollman (2014, Example~7.1.3) states the rules of the Fibonacci ambiguously.  In terms of Fibonacci numbers, the key rule is as follows: After winning a bet of $F_n$, stop if $n=1$ or 2, bet $F_{n-2}$ if $n\ge3$; after losing a bet of $F_n$, bet $F_{n+1}$.

Ethier (2010) and Bollman (2014) discuss the martingale, Fibonacci, Labou\-chere, and d'Alembert systems.  Hannum \& Cabot (2005), Epstein (2013), and Rodr\'iguez \& Mendes (2018) neglect only the Fibonacci, Packel (2006) discusses only the martingale and Labou\-chere systems, and Taylor (2015) discusses only the martingale, among these four betting systems.

\subsection{Gambler's ruin formula}

The gambler's ruin formula (for the probability of winning $W$ units before losing $L$ units when betting one unit at even money at each trial in an independent sequence) is important and easy to derive.  One can solve a linear system of second-order difference equations, the method used by Packel (2006, Section~7.4) and Taylor (2015, pp.~237--245).  Here one can point to the parallel with \textit{differential} equations.  If this seems too technical, one can use the martingale proof of De Moivre (Ethier, 2010, pp.~271--272), which is easy to understand even if not entirely rigorous.  Here is the argument.

\begin{quote}
Let us refer to the gambler as A and to his opponent as B.  We assume that A has $L$ units initially while B has $W$.  Play continues until one of the players has all $L+W$ units.  Let us call these units ``chips'' (De Moivre used the word ``counters'') and assume that each player stacks his chips and that each coup amounts to a transfer of one chip from the top of one stack to the top of the other.  The key idea is to artificially assign values to the chips that make the
game fair.  Thus, we assign values $q/p, (q/p)^2,\ldots,(q/p)^L$ to A's initial stack of chips from bottom to top, and $(q/p)^{L+1},\ldots,(q/p)^{L+W}$ to B's initial stack of chips from top to bottom.  At every coup, B risks $q/p$ times as much as A, and since A wins, loses, and pushes each coup with probabilities $p$, $q$, and $r:=1-p-q$, respectively, the game is fair.  Consequently, at the end of play, the expected (artificial) value of A's profit is
\begin{align*}
&\P(\text{A wins})\bigg\{\bigg(\frac{q}{p}\bigg)^{L+1}+\cdots+\bigg(\frac{q}{p}\bigg)^{L+W}\bigg\}\\
&\quad{}-(1-\P(\text{A wins}))\bigg\{\frac{q}{p}+\bigg(\frac{q}{p}\bigg)^2+\cdots+\bigg(\frac{q}{p}\bigg)^L\bigg\},
\end{align*}
which must be 0.
\end{quote}
The formula for $\P(\text{A wins})$ follows, with the cases $p\ne q$ and $p=q$ treated separately.  (A rigorous treatment requires the optional stopping theorem from martingale theory.)

Hannum \& Cabot (2005, p.~198) and Epstein (2013, pp.~58--59) state the gambler's ruin formula without proof.

\subsection{Other topics}

There are several other mathematical topics, not already covered in discrete probability, that should not be neglected in a course on the mathematics of gambling.

The Kelly system (Kelly, 1956) is taken up by only three of the seven textbooks, perhaps because its minor use of calculus puts it beyond the scope of most of the pre-calculus textbooks.  Just as the central limit theorem can be (and often is) included in pre-calculus statistics textbooks, the Kelly system can be discussed in even the most elementary course.  Indeed, Taylor (2015, pp.~230--237) shows how to do it.  Epstein (2013, pp.~61--62) gives brief coverage, and Ethier (2010, Chap.~10) provides more-thorough coverage.

Bold play is a beautiful topic, which can be discussed at an elementary level.  Proofs are a challenge, but providing motivation for the principal results is quite easy, especially if the gambler's ruin formula is available.  Bollman (2014, Section~7.6) gives excellent motivation but stops short of stating a theorem.  Epstein (2013, p.~60) states a theorem, though not quite accurately:  The maximum boldness strategy is not necessarily optimal in the presence of a house limit.  Ethier (2010, Chap.~9) gives extended coverage.

Shuffling is a suitable topic for a course on the mathematics of gambling.  Hannum \& Cabot (2005, pp.~261--264, 281--285) discusses shuffle tracking and statistical properties of shuffles; Epstein (2013, pp.~222--232) discusses periodicity of perfect shuffles and what he refers to as the amateur shuffle (a.k.a.\ the riffle shuffle), among other things; and Ethier (2010, Section 11.1) proves the theorem of Bayer \& Diaconis (1992) on riffle shuffles that led to their celebrated finding that ``seven shuffles suffice.''

Finally, martingales are covered by only one of the seven textbooks, perhaps because they depend on conditional expectation, which itself receives little coverage.  The one result from martingale theory that seems essential to a course on the mathematics of gambling is the result that says that no betting system applied to a sequence of fair or subfair wagers can result in a positive expectation.  Of course, some hypotheses are needed.  An elementary formulation of this result can be found in Thorp (1984, pp.~121--124).

\section{Mathematical level}

Here we compare the seven textbooks in terms of mathematical level.  Only three of the seven books state their prerequisites clearly.\medskip  

Packel (2006, p.~ix) writes, 
\begin{quote}
While the only formal mathematics background assumed is high school algebra, some enthusiasm for and facility with quantitative reasoning will also serve the reader well.  
\end{quote}

Ethier (2010, p.~v) writes, 
\begin{quote}
This is a monograph/textbook [$\dots$] intended for those already familiar with probability at the post-calculus, pre-measure-theory level. 
\end{quote}

Rodr\'iguez \& Mendes (2018, p.~xiii) write,
\begin{quote}
The material should be suitable for a college-level general education course for undergraduate college students who have taken an algebra or pre-algebra class.  In our experience, motivated high-school students who have taken an algebra course should also be capable of handling the material.
\end{quote}

With this information, we can rank the seven textbooks according to their mathematical level.  We start with the more advanced books.

\begin{itemize}
\item Ethier (2010) assumes some familiarity with basic probability, and is intended for an upper-division class in mathematics or statistics.  Proofs are given for most of the principal results.

\item Epstein (2013) covers few topics in depth, but there are a vast number of topics (some of which are unrelated to gambling).  Proofs, which appeared in the original edition (Epstein, 1967), were dropped in the revised edition (Epstein, 1977) and remain absent in the latest edition.  Although the book is largely concerned with discrete mathematics, it does use calculus more than just occasionally.  We conclude that the prerequisite for this textbook is calculus.  And some familiarity with basic probability is likely also necessary.
\end{itemize}

We believe that the other five textbooks could be successfully used in a pre-calculus course.  That does not mean, however, that they are all at the same mathematical level.

\begin{itemize}
\item Taylor (2015) is a user-friendly textbook that starts at the beginning and eventually gets into some rather sophisticated mathematics.  Examples of nonstandard topics include a combinatorial derivation of the seven-card poker-hand probabilities (where only the best five-card subset counts); probability problems that can be solved recursively;  Sicherman dice via generating functions; stationary distributions of finite Markov chains\footnote{\,It is incorrectly stated (p.~170) that these are always unique.  A necessary and sufficient condition for uniqueness of stationary distributions is uniqueness of recurrent classes.} with application to Monopoly and other board games; and an analysis of winning streaks in independent Bernoulli trials.  Although calculus is only rarely needed, this book seems to require of its readers greater mathematical maturity than the other four textbooks.

\item Bollman (2014) introduces many games not listed in Table \ref{gambling-content} (including 
punchboards, Die Rich, spider craps, twenty-six, Royal Roulette, Diamond Roulette, Riverboat Roulette, Double Action Roulette, crapless craps, barbooth, Double Dice, card craps, EZ Baccarat, Rupert's Island Draw, card slots, Multicolore, boule, Double Exposure, Super Fun 21, Spanish 21, Multiple Action 21, and Blackjack Switch) and analyzes most of them.  Only the more complicated topics (e.g., blackjack, video poker) lack derivations and rely on other sources for their conclusions.  Because of a greater emphasis on derivations, albeit elementary ones, we regard Bollman (2014) as more advanced than Hannum \& Cabot (2005).

\item Hannum \& Cabot (2005, p.~ix) clearly states that its primary intended audience is casino managers, or potential managers.  The prerequisites are similar to those of Packel (2006) and Rodr\'iguez \& Mendes (2018), but its level of mathematical sophistication puts it ahead of those two books on this list.  Compare, for example, their treatments of blackjack and of the central limit theorem.
\end{itemize}

The last two textbooks would be suitable for a general education class.  Both have a high school algebra prerequisite.

\begin{itemize}
\item Rodr\'iguez \& Mendes (2018) has much to recommend it.  It has perhaps the clearest analysis of the pass-line bet at craps of any of the seven textbooks.  It avoids proofs but often confirms conclusions with simulations in \textit{R}.  It studies game theory in some depth, permitting analysis of a poker-like game.

\item Packel's (2006) primary selling point is its broad accessibility.  The price of this accessibility is its sparse coverage of probability (omitting conditional probability, random variables, and variance, for example), but, like Rodr\'iguez \& Mendes (2018), it does provide an introduction to game theory with an application to a poker endgame situation.
\end{itemize}

\section{Conclusions}

Every book on discrete probability has many examples; indeed, it would not be possible to learn the subject without them.  Some of the simplest and least contrived examples arise in games of chance, perhaps not surprisingly because that is where the subject originated more than three centuries ago.

A course on the mathematics of gambling can be regarded as a course on discrete probability in which \textit{all} of the examples involve gambling.  Although some instructors may prefer a broader emphasis, to include games that are not ordinarily considered gambling games (e.g., Yahtzee), or to include games that are not even games of chance (e.g., checkers), there is no question that there exists abundant material even within the narrow confines of casino games.  And furthermore, there is complete flexibility with respect to the mathematical level of such a course, from high school level to graduate level, but usually at the university undergraduate level, lower division or upper division.

To justify our claim that examples from games of chance are less contrived than ``real-world'' examples, we cite a problem from the excellent probability textbook by Ross (2010, Problem 3.51):

\begin{quote}
A worker has asked her supervisor for a letter of recommendation for a new job. She estimates that there is an 80 percent chance that she will get the job if she receives a strong recommendation, a 40 percent chance if she receives a moderately good recommendation, and a 10 percent chance if she receives a weak recommendation. She further estimates that the probabilities that the recommendation will be strong, moderate, and weak are 0.7, 0.2, and 0.1, respectively.

(a) How certain is she that she will receive the new job offer?

(b) Given that she does receive the offer, how likely should she feel that she received a strong recommendation? a moderate recommendation? a weak recommendation?

(c) Given that she does not receive the job offer, how likely should she feel that she received a strong recommendation? a moderate recommendation? a weak recommendation?
\end{quote}

\noindent For comparison, here is a problem from Ethier (2010, Problem 1.20), rewritten slightly to make it comparable to the Ross problem.

\begin{quote}
Recall the rules of the pass-line bet at craps as well as the distribution of the various dice totals, given
by $\pi_j:=(6-|j-7|)/36$, $j=2,3,4,\ldots,12$.

(a) Find the probability that the pass-line bet is won.

(b) Regarding the probabilities $\pi_j$ as the prior probabilities for the come-out roll, find the corresponding posterior probabilities, given that the pass-line bet is won.  

(c) Regarding the probabilities $\pi_j$ as the prior probabilities for the come-out roll, find the corresponding posterior probabilities, given that the pass-line bet is lost.
\end{quote}

\noindent Both problems have the same structure.  But in the first problem three conditional probabilities and three unconditional probabilities are assumed without any real basis, whereas in the second problem we assume only fair dice.  The first problem is contrived, the second one is not.

\subsection{Advice to potential instructors}

Our first advice to instructors preparing to teach a course on the mathematics of gambling is to read the literature of the field as thoroughly as time permits but to do so skeptically.  The point is that, especially since the dawn of the Internet, the barrier to publishing in this field has been very low, so there exist many works that are mathematically unreliable.  Barnhart (1988), referring specifically to the history of roulette, describes what he calls a ``scholarly disaster,'' a situation where assertions become accepted once they are misstated often enough.  

This can be prevented if instructors independently confirm any results they borrow from other sources.  This is usually easy to do, but there are exceptions (e.g., blackjack, video poker, Texas hold'em).  In the exceptional cases, it is important to recognize which sources are reliable, and this can only be learned through experience.

Whichever textbook is chosen by the instructor, we recommend supplementing it with additional material, to make the course fit the instructor's and the students' interests.  Just because a book gives inadequate coverage to a particular topic should not rule it out as a potential textbook for the course.  In addition, the instructor should consider bringing in outside speakers, if it is feasible.  For example, in the Stanford University course mentioned in the Introduction, one lecture was given by a professional poker player, another by the head statistician for the California State Lottery, and a third by the author of the textbook used for the course.

Moreover, even for students in a mathematics department, facility with computing, in addition to logical reasoning skills, is both relevant and important.  Many of the more interesting problems and calculations in applications of probability to gambling require numerical work or programming, and such a course may encourage students to learn a suitable programming language.

We also point out that, in the event that the course attracts exceptional students, there are accessible research problems that can challenge their abilities.  The University of Utah course mentioned in the Introduction resulted in a publication (Vanniasegaram, 2006) by a UC Berkeley undergraduate, and the course at Stanford University resulted in the solution of a 12-year-old open problem by two graduate students (Han \& Wang, 2019).

We also encourage instructors to include material that can captivate a student's imagination and pique his or her interest in the subject.  Let us mention several of our favorite examples.  
\begin{itemize}
\item When discussing Texas hold'em, play scenes from the 1998 film \textit{Rounders} in which Mike McDermott (Matt Damon) competes against Teddy ``KGB'' (John Malkovich) in heads-up play, or the scene from the 2006 film \textit{Casino Royale} in which James Bond (Daniel Craig) check-raises Le Chiffre (Mads Mikkelsen) to win a \$115 million pot.

\item When discussing chemin de fer, read Chapters 9--13 of Fleming's (1953) \textit{Casino Royale}, at the climax of which James Bond competes against Le Chiffre at a hand of baccarat/chemin de fer with 32 million francs at stake.  (This is every bit as dramatic as the hold'em game just cited.)

\item When discussing card counting at blackjack, play scenes from the 2008 film \textit{21} or read excerpts from the book on which the film is based, Mezrich's (2002) \textit{Bringing Down the House: The Inside Story of Six M.I.T. Students Who Took Vegas for Millions}.  Thorp's (1966b) \textit{Beat the Dealer} is also compelling, especially for historical reasons.

\item When discussing craps, read Scoblete's (2007, Part 4; 2010) dramatic story of the Captain's 148-roll hand.  As we mentioned, this record was subsequently eclipsed by Ms.~DeMauro, but no such eyewitness account of the latter incident exists.

\item When discussing the gambler's ruin formula, read Johnson's (1990) ``Tale of a whale: Mysterious gambler wins, loses millions,'' the story of a freeze-out match between a Japanese billionaire, who would later be murdered, and an Atlantic City casino mogul, who would later be President of the United States.  See also Crowley (2016) for a more recent account.

\item When discussing lotteries, read Arratia et al.'s (2015) ``Some people have all the luck'' about how certain Florida State Lottery players won so regularly that they must have been up to something.  Or read the story of a retired couple who exploited a loophole in the Michigan State Lottery to win millions (Fagone, 2018).

\item When discussing video poker, read Poulsen's (2014) ``Finding a video poker bug made these guys rich\,---\,Then Vegas made them pay.''  This article suggests that, even if one finds an exploitable weakness is a casino game, one cannot legally exploit it.  Indeed, the same conclusion apparently applies to Phil Ivey's edge-sorting scheme at baccarat (Hawkins, 2017), a mathematical analysis of which appears in Dalton \& Hannum (2016).

\item When discussing the Labouchere system, read Leigh's (1976) \textit{Thirteen Against the Bank}, the story of how 13 Englishmen and women beat the roulette wheels of the casino at Nice using the reverse Labouchere system.  This book is alleged to be a work of nonfiction.  A good project for the class would be to investigate its plausibility, perhaps using simulation.

\item When discussing the Kelly system, read Poundstone's (2005) \textit{Fortune's Formula: The Untold Story of the Scientific Betting System that Beat the Casinos and Wall Street}.  It brings together a remarkable cast of characters, from mobsters to Nobel laureates.

\item When discussing roulette, read excerpts from Dostoevsky's 1867 novel, \textit{The Gambler}.  Alternatively, read the exploits of a number of successful biased wheel players in Barnhart's (1992, Chapters 4--10) \textit{Beating the Wheel}.

\item In the event that pari-mutuel betting on horses is discussed, read Chellel's (2018) ``The gambler who cracked the horse-racing code,'' a profile of Bill Benter, the most successful horse bettor of all time. 

\item In the event that faro is discussed, read Pushkin's 1834 short story, ``The Queen of Spades.''
\end{itemize}
In fact, there is so much great literature concerned with gambling that any list such as this must necessarily be incomplete.

Instructors interested in introducing a course on the mathematics of gambling should keep in mind that there could be some resistance on the grounds that such a course would encourage gambling.  As Emerson et al.\ (2009) pointed out, 
\begin{quote}
Researchers estimate that 3\% to 11\% of college students in the U.S. have a serious gambling problem that can result in psychological difficulties, unmanageable debt, and failing grades.
\end{quote}
Without attempting to minimize the problem, we believe that such objections lack merit.  A student who understands the law of large numbers, the gambler's ruin formula, the futility of betting systems, and the house advantages of the standard casino games is less likely to engage in games of pure chance (roulette, craps, keno, baccarat) than a student who lacks this knowledge.  A student who understands card counting at blackjack, hand-rank tables at video poker, and the principles of Texas hold'em, is perhaps more likely to participate in these skill-based games of chance than a student who lacks this knowledge, but to do so possibly as an advantage player, or in any case in a way unlikely to cause harm.  We are not suggesting that learning to gamble intelligently is a solution to problem gambling, only that it is unlikely to create a new problem gambler.

\section*{Acknowledgments}

The authors are grateful to Persi Diaconis for valuable comments on a preliminary draft, and to two reviewers for their helpful suggestions.  SNE was partially supported by a grant from the Simons Foundation (429675), and FMH was partially supported by a Discovery grant from NSERC.

\begin{newreferences}

\item Arratia, R., Garibaldi, S., Mower, L., \& Stark, P. B. (2015).  Some people have all the luck. \textit{Math.\ Magazine}, \textit{88} (3), 196--211.

\item Bachelier, L. (1914). \textit{Le jeu, la chance et le hasard}. Paris: Ernest Flammarion. Reprinted by \'Editions Jacques Gabay, Sceaux, 1993.

\item B\v arboianu, C. (2013). \textit{Probability Guide to Gambling: The Mathematics of Dice, Slots, Roulette, Baccarat, Blackjack, Poker, Lottery and Sport Bets}, Second Edition. Craiova, Romania: Infarom Publishing.

\item Barnhart, R. T. (1988). The invention of roulette.  In W.~R. Eadington (Ed.) \textit{Gambling Research: Proceedings of the Seventh International Conference on Gambling and Risk Taking, Vol.\ 1. Public Policy and Commercial Gaming Industries throughout the World} (pp.~295--330). Reno: Bureau of Business and Economic Research, University of Nevada.

\item Barnhart, R. T. (1992). \textit{Beating the Wheel}.  New York: A Lyle Stuart Book, Carol Publishing Group. 

\item Bayer, D. \& Diaconis, P. (1992). Trailing the dovetail shuffle to its lair.  \textit{Ann. Appl. Probab.}, \textit{2} (2), 294--313.

\item Bewersdorff, J. (2005). \textit{Luck, Logic, and White Lies: The Mathematics of Games}. Wellesley, MA: A. K. Peters, Ltd. 

\item Bollman, M. (2014). \textit{Basic Gambling Mathematics: The Numbers Behind the Neon}.  Boca Raton, FL: CRC Press, an imprint of Taylor \& Francis Group.  

\item Brown, B. H. (1919). Probabilities in the game of ``shooting craps''. \textit{Amer.\ Math.\ Monthly}, \textit{26} (8), 351--352.

\item Chellel, K. (2018). The gambler who cracked the horse-racing code.  \textit{Bloomberg Businessweek}, May 3.
\url{https://getpocket.com/explore/item/the-gambler-who-cracked-the-horse-racing-code}.

\item Chen, B. \& Ankenman, J. (2006). \textit{The Mathematics of Poker.} Pittsburgh: ConJelCo LLC.

\item Croucher, J. S. (2003). \textit{Gambling and Sport: A Statistical Approach}.  North Ryde, NSW: Macquarie University Lighthouse Press.

\item Croucher, J. S. (2006). Teach statistics? You bet!  \textit{Significance}, \textit{3} (1), 46--48.

\item Crowley, M. (2016). The whale that nearly drowned the Donald: How Trump schemed to win back millions from a high-rolling\,---\,and doomed\,---\,Japanese gambler.  \textit{Politico Magazine}, Feb.\ 14.  Retrieved October 26, 2019, from \url{https://www.politico.com/magazine/story/2016/02/japanese-gambler-donald-trump-213635}.

\item Dalton, T. \& Hannum, R. (2016). The mathematics of baccarat edge sorting. \textit{J. Gambling Bus.\ Econ.}, \textit{10} (3), 71--82.

\item de Finetti, B. (1937).  La pr\'evision:  Ses lois logiques, ses sources subjectives. \textit{Ann.\ Inst.\ H. Poincar\'e}, \textit{7} (1), 1--68. 
English translation: Kyberg Jr., H. E. (1992). Foresight: Its logical laws, its subjective sources. In  S. Kotz \& N. L. Johnson (Eds.) \textit{Breakthroughs in Statistics} (pp.~134--174).  New York: Springer.

\item De Moivre, A. (1738). \textit{The Doctrine of Chances: or, A Method of Calculating the Probabilities of Events in Play}, Second Edition.  London: H. Woodfall.

\item Dubins, L. E. \& Savage, L. J. (2014).  \textit{How to Gamble If You Must: Inequalities for Stochastic Processes}, Second Dover Ed.  Mineola, NY: Dover Publications, Inc.

\item Emerson, P. V. D. et al. (2009).  Addressing College Gambling: Recommendations for Science-Based Policies and Programs.  A Report of the Task Force on College Gambling Policies.  A Project of the Division on Addictions, The Cambridge Health Alliance, A Teaching Affiliate of Harvard Medical School, and the National Center for Responsible Gaming.  Retrieved October 26, 2019, from \url{https://www.divisiononaddiction.org/html/publications/College_Report_Full.pdf}.

\item Epstein, R. A. (1967). \textit{The Theory of Gambling and Statistical Logic}. New York: Academic Press.

\item Epstein, R. A. (1977). \textit{The Theory of Gambling and Statistical Logic}, Revised Edition. New York: Academic Press.

\item Epstein, R. A. (2013). \textit{The Theory of Gambling and Statistical Logic}, Second Edition.  Waltham, MA: Academic Press, an imprint of Elsevier.

\item Ethier, S. N. (2010). \textit{The Doctrine of Chances: Probabilistic Aspects of Gambling}.  Berlin--Heidelberg: Springer-Verlag.

\item Ethier, S. N. \& Gamez, C. (2013). A game-theoretic analysis of \textit{baccara chemin de fer}. \textit{Games}, \textit{4} (4), 711--737.

\item Ethier, S. N. \& Hoppe, F. M. (2010). A world record in Atlantic City and the length of the shooter's hand at craps.  \textit{Math.\ Intelligencer}, \textit{32} (4), 44--48.

\item Ethier, S. N., Kim, J. J., \& Lee, J. (2019). Optimal conditional expectation at the video poker game Jacks or Better.  \textit{UNLV Gaming Res.\ \& Rev. J.}, \textit{23} (1), 1--18.

\item Ethier, S. N. \& Lee, J. (2015a). On the three-person game \textit{baccara banque}.  \textit{Games}, \textit{6} (2), 57--78.

\item Ethier, S. N. \& Lee, J. (2015b). The evolution of the game of baccarat.  \textit{J. Gambling Bus.\ Econ.}, \textit{9} (2), 1--13.

\item Ethier, S. N. \& Lee, J. (2019). Snackjack:  A toy model of blackjack.  Retrieved October 26, 2019, from \url{https://arxiv.org/abs/1906.01220}.

\item Fagone, J. (2018).  Jerry and Marge go large: Gaming the lottery seemed as good a retirement plan as any.  \textit{Huffington Post}, Mar.~1.  Retrieved October 26, 2019, from \url{https://highline.huffingtonpost.com/articles/en/lotto-winners/}.

\item Ferguson, C. \& Ferguson, T. (2007). The endgame in poker.  In S. N. Ethier \& W. R. Eadington (Eds.) \textit{Optimal Play: Mathematical Studies of Games and Gambling} (pp.~79--106).  Reno: Institute for the Study of Gambling and Commercial Gaming.

\item Fleming, I. (1953). \textit{Casino Royale}. London: Jonathan Cape Ltd. 

\item Gibbs, R. A. \& Johnson, L. S. (1977). Developing college courses on gambling.  Unpublished manuscript (8~pp.), available at UNLV Lied Library Special Collections,  LB 2365 .G3 G52.

\item Gould, R. J. (2016). \textit{Mathematics in Games, Sports, and Gambling: The Games People Play}, Second Edition. Boca Raton, FL: CRC Press, an imprint of Taylor \& Francis Group.  

\item Griffin, P. A. (1999). \textit{The Theory of Blackjack: The Compleat Card Counter's Guide to the Casino Game of 21}, Sixth Edition.  Las Vegas: Huntington Press.

\item Guerrera, T. (2007). \textit{Killer Poker by the Numbers: The Mathematical Edge for Winning Play}. New York: Lyle Stuart, Kensington Publishing Corp.

\item Haigh, J. (2003). \textit{Taking Chances: Winning with Probability}, Second Edition. Oxford, UK: Oxford University Press.  

\item Han, Y. \& Wang, G. (2019).  Expectation of the largest bet size in the Labouchere system. \textit{Electron.\ Commun.\ Probab.}, \textit{24} (11), 1--10.

\item Hannum, R. (2007). The partager rule at roulette: Analysis and case of a million-euro win. In S. N. Ethier \& W. R. Eadington (Eds.) \textit{Optimal Play: Mathematical Studies of Games and Gambling} (pp.~325--344).  Reno: Institute for the Study of Gambling and Commercial Gaming.

\item Hannum, R. C. \& Cabot, A. N. (2001). \textit{Practical Casino Math}. Reno: Institute for the Study of Gambling and Commercial Gaming.

\item Hannum, R. C. \& Cabot, A. N. (2005). \textit{Practical Casino Math}, Second Edition. Reno: Institute for the Study of Gambling and Commercial Gaming; Las Vegas: Trace Publications.

\item Hannum, R. C. \& Cabot, A. N. (2012). \textit{Practical Casino Math}, Chinese translation. Taipei: Yang--Chih Book Co.

\item Hannum, R., Rutherford, M., \& Dalton, T. (2012). Economics of poker: The effect of systemic chance. \textit{J. Gambling Bus.\ Econ.}, \textit{6} (1), 25--48.

\item Hawkins, D. (2017).  What is `edge-sorting' and why did it cost a poker star \$10 million in winnings?  \textit{Washington Post}, Oct.~26.  Retrieved October 26, 2019, from \url{https://www.washingtonpost.com/news/morning-mix/wp/2017/10/26/what-is-edge-sorting-and-why-did-it-cost-a-poker-star-10-million-in-winnings/}.

\item Henze, N. \& Riedwyl, H. (1998). \textit{How to Win More: Strategies for Increasing a Lottery Win}.  Natick, MA: A. K. Peters, Ltd.

\item Johnson, L. S. (1975). Las Vegas: Springboard to mathematical study\,---\,Fact or fiction?  Unpublished manuscript (10~pp.), available at UNLV Lied Library Special Collections,  GV 1301 .J5.

\item Johnson, L. S. (1977). Las Vegas\,---\,A classroom for colleges.  Unpublished manuscript (7~pp.), available at UNLV Lied Library Special Collections,  GV 1301 .J52.

\item Johnson, R. (1990). Tale of a whale: Mysterious gambler wins, loses millions.  \textit{Wall Street J.}, June 28, A1.

\item Kelly, J. L., Jr. (1956). A new interpretation of information rate. \textit{Bell System Tech.\ J.}, \textit{35} (4), 917--926.

\item Kemeny, J. G. \& Snell, J. L. (1957). Game-theoretic solution of baccarat.  \textit{Amer.\ Math.\ Monthly}, \textit{64} (7), 465--469. 

\item Laplace, P. S. (1819).  \textit{Th\'eorie analytique des probabilit\'es, Vol.~1}. Essai philosophique sur les probabilit\'es, 4$^e$ \'ed.  Paris: M. V. Courcier, Imprimeur-Libraire pour les Math\'ematiques.

\item Laplace, P. S. (1902).  \textit{A Philosophical Essay on Probabilities}, translated from the Sixth French Edition by F. W. Truscott.  London: Chapman \& Hall, Ltd.

\item Leigh, N. (1976). \textit{Thirteen Against the Bank}. New York: William Morrow and Co.

\item Livingston, A. D. (1968). `Hold me': a wild new poker game $\ldots$ and how to tame it. \textit{Life}, Aug.~16, 38--42.

\item Maitra, A. P. \& Sudderth, W. D. (1996). \textit{Discrete Gambling and Stochastic Games.} New York: Springer-Verlag. 

\item Manson, A. R., Barr, A. J., \& Goodnight, J. H. (1975). Optimum zero-memory strategy and exact probabilities for 4-deck blackjack. \textit{Amer.\ Statist.}, \textit{29} (2), 84--88.  Correction \textit{29} (4), 175.

\item Mezrich, B. (2002). \textit{Bringing Down the House: The Inside Story of Six M.I.T. Students Who Took Vegas for Millions}. New York: Free Press,  Simon \& Schuster.

\item Montmort, P. R. (1708). \textit{Essay d'analyse sur les jeux de hazard}. Paris: Chez Jacque Quillau. Published anonymously. 

\item Montmort, P. R. (1713). \textit{Essay d'analyse sur les jeux de hazard}, Second Edition. Paris: Chez Jacque Quillau. Published anonymously. 

\item Packel, E. W. (1981). \textit{The Mathematics of Games and Gambling}.  Washington, D.C.: The Mathematical Association of America.

\item Packel, E. W. (2006). \textit{The Mathematics of Games and Gambling}, Second Edition.  Washington, D.C.: The Mathematical Association of America.

\item Poulsen, K. (2014).  Finding a video poker bug made these guys rich\,---\,Then Vegas made them pay. \textit{Wired}, Oct.~7.  Retrieved October 26, 2019, from \url{https://www.wired.com/2014/10/cheating-video-poker/}.

\item Poundstone, W. (2005). \textit{Fortune's Formula: The Untold Story of the Scientific Betting System that Beat the Casinos and Wall Street}. New York: Hill and Wang, a division of Farrar, Straus, and Giroux.

\item R\'enyi, A. (1972). \textit{Letters on Probability}. Detroit: Wayne State University Press.

\item Rodr\'iguez, A. \& Mendes, B. (2018). \textit{Probability, Decisions and Games: A Gentle Introduction Using R}.  Hoboken, NJ: John Wiley \& Sons.

\item Ross, S. (2010). \textit{A First Course in Probability}, Eighth Edition.  Upper Saddle River, NJ: Pearson Prentice Hall.

\item Schneider, W. \& Turmel, J. (1975). On the organization of a university level course in gambling.  Unpublished manuscript (15~pp.), available at UNLV Lied Library Special Collections, LB 2365 .M3 S35.  

\item Scoblete, F. (2007). \textit{The Virgin Kiss and Other Adventures.} Daphne, AL: Research Services Unlimited.  

\item Scoblete, F. (2010). The Captain rolls 147.  \textit{Golden Touch Craps}.  Retrieved October 26, 2019, from \url{https://www.goldentouchcraps.com/captainrolls147.shtml}.

\item Selvin, S. (1975a). A problem in probability. \textit{Amer.\ Statistician}, \textit{29} (1), 67.

\item Selvin, S. (1975b). On the Monty Hall problem. \textit{Amer.\ Statistician}, \textit{29} (3), 134.

\item Shackleford, M. (2016). Texas hold'em\,---\,Top hands for 10-player game.  \textit{Wizard of Odds}.  Retrieved October 26, 2019, from \url{https://wizardofodds.com/games/texas-hold-em/10-player-game/}.

\item Shackleford, M. (2019a). \textit{Gambling 102: The Best Strategies for All Casino Games}, Second Edition. Las Vegas: Huntington Press.

\item Shackleford, M. (2019b). Craps side bets.  \textit{Wizard of Odds}.  Retrieved October 26, 2019, from \url{https://wizardofodds.com/games/craps/appendix/5/}.

\item Suddath, C. (2009). Holy craps! How a gambling grandma broke the record. \textit{Time}, May~29.  Retrieved October 26, 2019, from \url{http://content.time.com/time/nation/article/0,8599,1901663,00.html}.

\item Taylor, D. G. (2015). \textit{The Mathematics of Games: An Introduction to Probability}. Boca Raton, FL: CRC Press, an imprint of Taylor \& Francis Group.  

\item Thorp, E. O. (1966a).  \textit{Elementary Probability}. New York: John Wiley \& Sons.

\item Thorp, E. O. (1966b). \textit{Beat the Dealer: A Winning Strategy for the Game of Twenty-One}, Revised and Simplified Edition. New York: Random House.

\item Thorp, E. O. (1984). \textit{The Mathematics of Gambling}.  Secaucus, NJ: A Gambling Times Book, Lyle Stuart.

\item Vanniasegaram, S. (2006). Le her with $s$ suits and $d$ denominations. \textit{J. Appl.\ Probab.}, \textit{43} (1), 1--15.

\item Werthamer, N. R. (2018). \textit{Risk and Reward: The Science of Casino Blackjack}, Second Edition.  Cham, Switzerland: Springer International Publishing AG.

\item Wilson, A. N. (1970). \textit{The Casino Gambler's Guide}, Extended Edition.  New York: Harper \& Row.

\end{newreferences}

\end{document}